\documentclass[11pt]{amsart}

\usepackage{amsthm, amsfonts, amssymb, amsmath, latexsym, enumerate}
\usepackage[all]{xy}
\usepackage{mathrsfs}
\usepackage{array} 
\usepackage{hyperref}

\newtheorem{thm}{Theorem}[section] \newtheorem{lem}[thm]{Lemma}
\newtheorem{corol}[thm]{Corollary} \newtheorem{prop}[thm]{Proposition}

\newtheorem*{thm*}{Theorem}

\theoremstyle{definition} \newtheorem{rmk}[thm]{Remark}
 \newtheorem{dfn}[thm]{Definition}
 
\newtheorem{claim}[thm]{Claim}
 
 \newtheorem*{akn}{Acknowledgments}

\newcommand{\OO}{\mathscr{O}} \newcommand{\FF}{\mathscr{F}}
\newcommand{\GG}{\mathscr{G}} \newcommand{\EE}{\mathscr{E}}
 
\newcommand{\KK}{\mathscr{K}} \newcommand{\CC}{\mathscr{C}}
 \newcommand{\DD}{\mathscr{D}}
\newcommand{\EEnd}{\mathscr{E}nd} \newcommand{\EExt}{\mathscr{E}xt}

\newcommand{\Most}{{\sf M^s}}
\newcommand{\Moss}{{\sf M^{s\,s}}}

\newcommand{\FMst}{{\sf FM^{s}}}

\newcommand{\MCMst}{{\sf MCM^s}}
\newcommand{\FMCMst}{{\sf FMCM^s}}
\newcommand{\MCMss}{{\sf MCM^{s\,s}}}

\DeclareMathOperator{\VV}{{\mathbb V}} 
\DeclareMathOperator{\rk}{rk} 
\DeclareMathOperator{\coker}{coker} \DeclareMathOperator{\Hilb}{Hilb}

\DeclareMathOperator{\ts}{\otimes} \DeclareMathOperator{\SL}{{\sf SL}}
\DeclareMathOperator{\GL}{{\sf GL}}
 \DeclareMathOperator{\Pic}{Pic}

 \DeclareMathOperator{\lequiv}{=}

\newcommand{\Z}{\mathbb Z} \newcommand{\C}{\mathbb C}
 
 \newcommand{\p}{\mathbb P}
 \newcommand{\G}{\mathbb G}

\DeclareMathOperator{\Ext}{Ext} 
 
\DeclareMathOperator{\Hom}{Hom} \DeclareMathOperator{\im}{Im}
\DeclareMathOperator{\cok}{coker} \DeclareMathOperator{\End}{End}
\DeclareMathOperator{\HH}{H} \DeclareMathOperator{\hh}{h}
\DeclareMathOperator{\ext}{ext} 

\DeclareMathOperator{\len}{len}
\newcommand{\one}{\mathrm{id}}

\newcommand{\rr}{\rightarrow} 

 \newcommand{\xr}{\xrightarrow}
\newcommand{\LL}{\mathscr{L}} 
\newcommand{\MM}{\mathscr{M}} \newcommand{\NN}{\mathscr{N}}
\newcommand{\gnrt}{{\sf Gen}} \newcommand{\gen}{{\sf g}} \newcommand{\gnrtrs}{{\sf G}}
\newcommand{\syz}{{\sf Syz}}
\newcommand{\kk}{\C}
\newcommand{\fff}{{\sf f}} \newcommand{\ppp}{{\sf p}}

\newcommand{\LX}{{\sf L}(X)}
\newcommand{\CX}{{\sf C}(X)}
\newcommand{\TX}{{\sf T}(X)}
\newcommand{\RdX}{{\sf R}_d(X)}

\SelectTips{cm}{} \newdir{ >}{{}*!/-5pt/@{>}}
\numberwithin{equation}{section}
\allowdisplaybreaks[4]

\begin{document}


\title[aCM bundles on the cubic surface]{Rank 2 arithmetically Cohen-Macaulay bundles on a nonsingular cubic surface}

\author{Daniele Faenzi} \email{{\tt faenzi@math.unifi.it}}
\address{Dipartimento di Matematica ``U.~Dini'', Universit\`a di
Firenze, Viale Morgagni 67/a, I-50134, Florence, Italy}
\urladdr{\tt {\url{http://www.math.unifi.it/~faenzi/}}}


\thanks{The author was partially supported by Italian MIUR funds.}




\sloppy

\begin{abstract}
Rank $2$ indecomposable arithmetically Cohen-Macaulay bundles $\EE$ on
a nonsingular cubic 
surface $X$ in $\p^3$ are classified, by means of the possible forms taken by
the minimal graded free resolution of $\EE$ over $\p^3$.
The admissible values of the Chern classes of $\EE$ are listed and the
vanishing locus of a general section of $\EE$ is studied. 

Properties of $\EE$ such as slope (semi) stability and simplicity are
investigated; the number of relevant families is computed together
with their dimension.
\end{abstract}

\maketitle



\section{Introduction}

Given a smooth projective variety $Y$ of dimension $n$, equipped with a very ample line bundle
$\OO_Y(1)$, a vector bundle $\EE$ on $Y$ is called arithmetically
Cohen-Macaulay (aCM) if all its intermediate cohomology modules
vanish, i.e. if $\HH^p(Y,\EE(t))=0$ for $p\neq 0,n$ and for all $t \in
\Z$.

The set of aCM bundles on projective varieties has been studied in a large number of papers.
The splitting criterion of Horrocks, cfr. \cite{horrocks:punctured},
asserts if $Y$ is a projective space, then $\EE$ splits as a sum of line bundles.
Kn\"orrer in \cite{knorrer:aCM} proved that if $Y$ is a smooth quadric, then $\EE$
is a direct sum of line bundles and twisted spinor bundles.
The connection of these {\em splitting criteria} with the structure of the derived category
has been explored in \cite{ancona-ottaviani:beilinson-quadrics}.
The link with liaison theory should also be mentioned, cfr. the papers
\cite{hartshorne-casanellas:liaison}, 
\cite{hartshorne-drozd-casanellas:liaison}.

If there exists on $Y$ a 
finite set of isomorphism classes of aCM indecomposable bundles (up to twist by $\OO_Y(t)$)
then $Y$ is called of {\em finite Cohen-Macaulay type}.
It turns out that these varieties are completely classified, cfr. \cite{eisenbud-herzog} and reference therein.

The question was then posed of studying families of aCM bundles, at least those of low rank (say rank $2$), on 
varieties which are not of finite Cohen-Macaulay type.
The majority of results in this direction starts from the assumption that $\Pic(Y)\simeq \Z$.
For instance, the case of prime Fano threefolds has been analyzed in \cite{arrondo-costa},
\cite{madonna:fano-cy}, \cite{V5}, \cite{v22:pjm-to-appear}; see also \cite{arrondo-grana:G14},
for the case of $\G(\p^1,\p^4)$. 
In a similar spirit, general hypersurfaces of dimension $n \geq 3$ 
have been studied in \cite{chiantini-madonna:sextic}, \cite{chiantini-madonna:general}, \cite{kumar-rao-ravindra:hypersurfaces}.

As it results from \cite{arrondo-costa},
the general cubic threefold $Y$ admits three families of indecomposable aCM rank $2$ bundles,
corresponding respectively to a line, a conic, an elliptic quintic in $Y$.
Cutting with a hyperplane gives $3$ families of aCM bundles on a cubic surface $X$.
But are these the only families? What is their dimension?
Do at least the Chern classes of an arbitrary indecomposable rank $2$ aCM bundle $\EE$ lift to $\p^3$?
The aim of this paper is to classify completely these bundles, whereby
answering to many questions of this sort. 
The problem gets more intricate due to the rich structure of $\Pic(X)$.
We find $12$ types of bundles, i.e. $9$ more than the ones mentioned above.
For every single type we study the relevant families in terms of moduli spaces.
A brief summary of our results is the following:

\begin{thm*}
If $\EE$ is a rank $2$ indecomposable arithmetically Cohen-Macaulay bundle on a
nonsingular cubic surface $X$, then $\EE$ is one of the following $12$ types:
\begin{scriptsize}
\begin{equation} \label{all}
\begin{array}{c|c|c}
\mbox{\em Minimal Free Resolution} & \mbox{\em Chern} & \mbox{\em Families} \\
\hline
\hline
\begin{array}{c|c|c}
\mbox{Ref.} & \gnrt(\EE)  & \syz(\EE) \\
\hline
\hline
(\ref{6}.\ref{6-1}) & \OO^6 & \OO(-1)^6 \\
(\ref{6}.\ref{6-2}) & \OO^6 & \OO(-1)^6 \\
(\ref{6}.\ref{6-3}) & \OO^6 & \OO(-1)^6 \\
\hline
(\ref{5-2}.\ref{5-1-1}) & \OO^5 & \OO(-1)^4  \oplus \OO(-2) \\
(\ref{5-2}.\ref{5-1-2}) & \OO^5 & \OO(-1)^4  \oplus \OO(-2) \\
(\ref{5-4}.\ref{5-2-1}) & \OO \oplus \OO(-1)^4 & \OO(-2)^5  \\
(\ref{5-4}.\ref{5-2-2}) & \OO \oplus \OO(-1)^4 & \OO(-2)^5  \\
\hline
(\ref{4-1}) & \OO^4 & \OO(-1)^2 \oplus \OO(-2)^2 \\
(\ref{4-3}) & \OO^3 \oplus \OO(-1) & \OO(-1) \oplus \OO(-2)^3 \\
(\ref{4-2}) & \OO^2 \oplus \OO(-1)^2 & \OO(-2)^4 \\
(\ref{4-4}) & \OO \oplus \OO(-1)^3 & \OO(-2)^3 \oplus \OO(-3) \\
\hline
(\ref{3}) & \OO^3 & \OO(-2)^3 
\end{array}
& 
\begin{array}{c|c}
c_1 & c_2 \\
\hline
\hline 
2\,H & 5 \\
H+T & 4 \\
H+C+L & 3 \\
\hline
H+C & 3 \\
H+L_1+L_2 & 2 \\
H-C & 1 \\
H-L_1-L_2 & 0 \\
\hline
H+L & 2 \\
H & 2 \\
C & 1 \\
0 & 1 \\
\hline
T & 1 \\
\end{array}
& 
\begin{array}{c|c|c}
\mbox{\em num} & \mbox{\em dim} & \mbox{\em stab} \\
\hline
\hline
1 & 5 & ss - st\\
72 & 3 & ss - st \\
270 & 1 &ss - st\\
\hline
27 & 2 & u - st\\
216 & 0 & st \\
27 & 2 & u - st\\
216 & 0 & st\\
\hline
27 & 1 & ss - st\\
1 & 2  & u - st \\
27 & 1 & ss - st\\
1 & 2  & u - ss\\
\hline
72 & 0 & st
\end{array}
\end{array}
\end{equation}
\end{scriptsize}
\end{thm*}

Here $\gnrt(\EE)$ (resp. $\syz(\EE)$) describes the set of generators (resp. syzygies)
in the minimal graded
free resolution of (the extension by zero to $\p^3$ of)
$\EE$.
{\em Chern} outlines the Chern classes of $\EE$, where $L$, $L_i$, $C$, $T$ are
divisor classes corresponding respectively to lines, conics, twisted cubics contained in $X$.
{\em Num.}, {\em dim.}, indicate the
number respectively the dimension of each family.
The column {\em stab.} tells whether we can find an unstable (u),
strictly semistable (ss), and a stable (st) bundle $\EE$ with the prescribed invariants.
We will see that unstable bundles correspond to a finite number of instances.

The paper is organized as follows. In the next section we set up some background.
In section \ref{line-bundles}, we start by classifying aCM line
bundles on $X$.
In sections \ref{clas-res} and \ref{chern} we analyze the form of the
minimal graded free resolution of the aCM vector bundle $\EE$ and its
Chern classes.

Of course, if the bundle $\EE$ is an extension of two
aCM line bundles, it must be comprehended by our analysis.
We will focus on this
in section \ref{ext}, where
we prove that all cases of table (\ref{all}) contain indecomposable 
extension bundles, also called {\em layered sheaves}, cfr. \cite{hartshorne-drozd-casanellas:liaison}.

However not all of the aCM bundles are extensions: we will clarify this 
in sections \ref{moduli} and \ref{numerology}, where we study the
moduli spaces and the number of distinct families of aCM bundles of
rank $2$.

The material we need concerning the nonsingular cubic surface 
is contained in 
\cite{hartshorne:ag}, \cite{griffiths-harris:ag},
\cite{segre:cubic}, \cite{manin:cubic-russian}.
For the reader's convenience, we recall in the Appendix some basic facts about the
combinatorics of divisors classes on smooth cubic surfaces.

\begin{akn}
I would like to thank R. Hartshorne, whom I met in Torino at the
conference Syzygy 2005 in honor of P. Valabrega, for many valuable remarks,
and for pointing out to me the interesting papers
\cite{baciu-ene-pfister-popescu:rank-2} and \cite{ene-popescu:rank-1},
where similar questions are investigated.
Also I would like to thank G. Ottaviani for many inspiring ideas
and for his constant support.
\end{akn}

\section{Generalities} \label{generalities}

The material contained in this section is well-known.
We will work on a projective variety $Y$ over the field $\kk$ of complex numbers,
equipped with a very ample line bundle $\OO_Y(1)$, associated to the hyperplane class $H$.
For any $t \in \Z$, the line bundle
$\OO_Y(1)^{\otimes \,t}$ will be denoted indifferently by
$\OO_X(t)$ or $\OO_X(t \, H)$. 

Given a subvariety $Z$ of a smooth variety $Y$, we denote the ideal sheaf
(resp. the normal sheaf) of $Z$ in $Y$ by $J_{Z,Y}$ (resp. by $N_{Z,Y}$).
We will drop the subscript $Y$ whenever possible.
Given a positive integer $m$, we denote the Quot-scheme parameterizing 
subscheme $Z\subset Y$ of length $m$ by $\Hilb_{m}(Y)$
(see \cite[pag. 41]{huybrechts-lehn:moduli}).

Given a line bundle $\LL$, we write $|\LL|$ for the
linear system of sections of $\LL$, $\LL(D)$ for a twist of the line
bundle $\LL$ by the divisor $D$.  
If $\FF_1$, $\FF_2$ are coherent sheaves on $Y$, we will write $\hh^p(\FF_1)$,
$\ext^p(\FF_1,\FF_2)$, etc. for
the dimension over $\kk$ of the vector spaces $\HH^p(Y,\FF_1)$, $\Ext^p(Y;\FF_1,\FF_2)$.
According to Chiantini and Madonna, cfr. \cite{chiantini-madonna},
we will say that a torsionfree sheaf $\FF$ on $Y$ is {\em normalized} (with
respect to $\OO_Y(1)$) if
$\hh^0(Y,\FF(-t))=0$ for $t >0$, but $\hh^0(Y,\FF)\neq 0$.

From now on, $X$ will denote a {\em smooth cubic surface} in $\p^3$,
defined by a cubic form $F$, over the field $\kk$.
The very ample line bundle $\OO_X(1)$ will be the restriction of $\OO_{\p^3}(1)$,
and we will consider slope (semi) stability of sheaves on $X$
with respect to the hyperplane polarization $H$.
The Chern classes of a sheaf over $X$ will be considered as elements of $\HH^*(X,\Z)$.
In particular $c_1$ will belong to $\HH^{2}(X,\Z) \cong \Z^7$,
an abelian group equipped with an action of the Weil group ${\sf W}({\mathbf E}_6)$ (cfr. the Appendix),
while $c_2$ will be indicated by an integer.

We will denote the moduli space of rank $r$ slope-stable
(resp. semistable)
vector bundles $\EE$ on $X$ with
$c_1(\EE)=c_1$, $c_2(\EE)=2$ by $\Most(r;c_1,c_2)$
(resp. by $\Moss(r;c_1,c_2)$).
We refer the reader to \cite[Chapter 8.2]{huybrechts-lehn:moduli} for
the construction of this space.

A sheaf $\FF$ on $X$ is {\em simple} (resp. {\em rigid}, {\em unobstructed}) if $\Hom(\FF,\FF) = \kk$,
(resp. $\Ext^1(\FF,\FF) = 0$, $\Ext^2(\FF,\FF) = 0$).
Notice that a simple bundle is {\em indecomposable}.
Irreducibility of moduli spaces on Del Pezzo surfaces has
been analyzed by Gomez in \cite{gomez:thesis}.
A sheaf $\FF$ is called a {\em split} bundle if it decomposes as a direct sum of line bundles.

\begin{dfn}
A vector bundle $\EE$ on $X$ is aCM (i.e. arithmetically
Cohen-Macaulay) if it has {\em no intermediate cohomology}, i.e. if:
\begin{equation} \label{def-acm}
 \HH^1(X,\EE (t)) = 0 \qquad \mbox{for all $t \in \Z$}
\end{equation}

Of course this is an open condition.
So we denote the open subset of (a union of components of) the moduli space of stable 
(resp. semistable) vector bundles $\EE$ on $X$ of rank $r$ with
$c_1(\EE)=c_1$, $c_2(\EE)=c_2$ consisting of aCM sheaves by
$\MCMst_X(r;c_1,c_2)$ (resp. by $\MCMss_X(r;c_1,c_2)$).
\end{dfn}

We will consider also the moduli space $\FMst(r;c_1,c_2)$ of
{\em framed} stable sheaves i.e. pairs $[\EE,s]$ where $[\EE]$ is the class
of a sheaf in $\Most(r;c_1,c_2)$ and $[s]$ is an element of
$\p(\HH^0(\EE))$.
If $\HH^0(\EE)\neq 0$, we have a rational map
$\eta: \FMst(r;c_1,c_2) \rr \Most(r;c_1,c_2)$,
dominating any irreducible component containing $[\EE]$, with $\eta^{-1}([\EE])=\p(\HH^0(\EE))$.
Depending on the purpose, we will consider semistable, stable, aCM
framed sheaves, with obvious notation.

\subsection{Bundles on hypersurfaces}

The following theorem is well known, and a proof can be found in \cite{beauville:determinantal}.
If $\varepsilon \in \{-1,1\}$, we say that a matrix $f$ is $\varepsilon$-symmetric 
if $f^{\top}=\varepsilon \, f$.
Correspondingly we have a notion of {\em $\varepsilon$-symmetric duality} on
a vector bundle $\FF$, namely an isomorphism $\kappa:\FF \rr \FF^*(t)$
such that $\kappa^{\top} = \varepsilon \, \kappa$.

\begin{thm} \label{res-general}
Let $Y=\mathbb{V}(F_Y)$ be a smooth hypersurface of degree $d$ in $\p^n$ and let $\FF$ be an aCM
rank $r$ vector bundle on $Y$.
Then the minimal graded free resolution of the sheaf $\FF$, extended by zero to
$\p^n$, takes the form:
   $$0 \rr \syz(\FF) \xrightarrow{\fff(\FF)} \gnrt(\FF) \xr{\ppp(\FF)} \FF \rr 0$$
with $\gnrt(\FF) = \bigoplus\limits_{i=1}^{s} \OO_{\p^n}(b_i)$, 
$\syz(\FF) = \bigoplus\limits_{j=1}^{s} \OO_{\p^n}(a_j)$ and $\det(\fff(\FF))=F^r$.

Moreover, suppose that there exists an $\varepsilon$-symmetric
duality $\kappa: \FF \ts \FF \rr \OO_{Y}(d+t)$.
Then we have a natural
isomorphism $\syz(\FF) \simeq \gnrt(\FF)^*(t)$,
and $\fff(\FF) = \varepsilon \,\, \fff(\FF)^{\top}$.
\end{thm}

We will order the integers $a_i$'s and $b_j$'s so that $a_1\geq \cdots \geq a_s$,
$b_1\geq \cdots \geq b_s$. The $(i,j)$-th entry of the matrix $\fff(\FF)$ has degree
$b_j-a_i$, and we will sometimes write $\fff(\FF)$ as a matrix $(\alpha_{i,j})$ of
integers $\alpha_{i,j}=b_j-a_i$.

\begin{rmk} \label{periodicity} 
Let $Y \subset \p^n$ be as above, let $\deg(Y)=d$, and let $\LL$ (resp. $\FF$) be a line bundle
(resp. a rank 2 vector bundle) on $Y$.
Then Theorem \ref{res-general} implies:
\begin{enumerate}
\item The matrix $\fff(\LL)$ is symmetric iff $\exists \, t\in \Z$ with $\LL^{\ts 2} \simeq \OO_{Y}(t)$;
\item The matrix $\fff(\FF)$ is skew-symmetric iff $\exists \, t\in \Z$ with $\wedge^2(\FF) \simeq \OO_{Y}(t)$;
\item \label{minors} If $\rk(\gnrt(\FF))=s$, then any minor if order $(s-1)$ of $\fff(\FF)$ vanishes on $Y$
  (i.e. any such minor is divided by $F$). 
\end{enumerate}

Moreover, by a result of Eisenbud (cfr. \cite{eisenbud:homological}),
if $f:\gnrtrs_1 \rr \gnrtrs_0$ is a presentation matrix over $\p^n$
for the vector bundle $\FF$ on $Y$
(i.e. $\coker(f) \simeq \FF$), then
there exists an infinite $2$-periodic exact sequence (perhaps non minimal) of
the form:

\begin{equation} \label{periodic}
\cdots \rr
\gnrtrs_3 \ts \OO_Y \xr{f_{|Y}} \gnrtrs_2 \ts \OO_Y
\xr{g_{|Y}} \gnrtrs_1 \ts \OO_Y \xr{f_{|Y}} \gnrtrs_0 \ts \OO_Y
\rr \FF \rr 0
\end{equation}
with $\gnrtrs_{2\,k} = \gnrtrs_0(-k\,d)$,
$\gnrtrs_{2\,k + 1} = \gnrtrs_1(-k\,d)$,
and where the map $g: \gnrtrs_2 \rr \gnrtrs_1$
is a resolution matrix over $\p^n$ for $\ker(f_{|Y})$, which is
an aCM vector bundle on $Y$ of rank $\rk(\gnrtrs_0)-\rk(\FF)$.
Therefore $g$ gives a resolution of the syzygy bundle $\ker(\fff(\EE)_{|Y})$.
Notice that a resolution is necessarily minimal if there is no {\em constant morphism} (i.e. no map of degree $0$) between
any summands of $\syz(\EE)$ and $\gnrt(\EE)$.

One should also notice a converse to Theorem \ref{res-general}, namely given a square matrix $\fff$ on $\p^n$
between split bundles of rank $s$, with $\det(\fff)=F^r$, if:
\begin{enumerate}
\item all minors of order $s-r+1$ of $\fff$ vanish on $\VV(F)$,
\item at any point there is a nonzero minor of order $s-r$ of $\fff$,
\end{enumerate}
then $\cok(\fff)$ is a rank $r$ aCM bundle on $\VV(F)$.
\end{rmk}

\subsection{Codimension 2 subschemes}

The Serre correspondence relates rank $2$ vector bundles on $X$
to subschemes $Z\subset X$ of codimension $2$.

For the proof of the following theorem we refer to \cite[Theorem 5.1.1]{huybrechts-lehn:moduli}.

\begin{thm} \label{hartshorne-serre}
Let $Z\subset X$ be a locally complete intersection subscheme of codimension $2$ in $X$,
and let $\LL$ be a line bundle on $X$.
Then the following are equivalent:
\begin{enumerate}[i)]
\item \label{CB-1} There exist a vector bundle $\EE$ with $\wedge^2 \EE \simeq \LL$ and an extension:
\begin{equation} \label{CB}
0 \rr \LL^* \rr \EE^* \rr J_Z \rr 0
\end{equation}
\item \label{CB-2} The pair $(\LL\ts \omega_X,Z)$ has the Cayley-Bacharach property i.e.
for any $s\in \HH^0(\LL\ts \omega_X)$, and for any $Z'\subset Z$ with $\len(Z')=\len(Z)-1$,
we have $s_{|Z}=0 \Leftrightarrow s_{|Z'}=0$.
\end{enumerate}
\end{thm}

Notice that dualizing (\ref{CB}) we obtain the exact sequence:
\begin{equation} \label{CB-dual}
0 \rr \OO_X \xr{s} \EE \rr J_Z \ts \det(\EE) \rr 0
\end{equation}

We will make use of the following Remark.
The proof of the statements regarding $\Hilb_{m}(X)$ can be found e.g. in
\cite[pag. 104]{huybrechts-lehn:moduli}.

\begin{rmk} \label{hilbert-serre}
The vector bundle $\EE^*$ of the previous theorem provides an
extension class which is an element of $\Ext^1(J_Z,\LL^*)$.
By Serre duality we have:
\begin{equation}
\Ext^1(J_Z,\LL^*)^* \simeq \HH^1(J_Z\ts \LL \ts \omega_X)
\end{equation}

Set $\ell = \len(Z)$.
Whenever $\Ext^1(J_Z,\LL^*)\simeq \kk$, we associate to $Z\subset X$ a unique pair $(\EE_Z,s_Z)$, where
$\EE_Z$ fits in the extension (\ref{CB}), $s_Z\in \HH^0(\EE_Z)$ and $Z=\{s_Z=0\}$.
If $\EE_Z$ is stable, this defines a map locally around $[Z] \in \Hilb_{\ell}(X)$:
  $$\zeta: \Hilb_{\ell}(X) \dashrightarrow \FMst(2;c_1(\LL^*),\ell)$$

We also have a rational map defined around the point $[\EE_Z,s_Z]$, namely we associate to a section its zero locus:
  $$\xi: \FMst(2;c_1(\LL^*),\ell) \dashrightarrow \Hilb_{\ell}(X)$$

The map $\xi$ is dominant if and only if 
the pair $(\EE_Z,s_Z)$ is defined for a general subscheme
$Z$ of $X$ of length $\len(Z)$.
If $\EE_Z$ is a simple bundle (i.e. $\End(\EE) \simeq \C$),
then $\xi$ is birational onto its image, $\zeta$ being its local inverse.
\end{rmk}

\section{Line bundles} \label{line-bundles}

We will use the notation $\RdX$ for
the set of divisor classes containing degree $d$ smooth irreducible rational curves contained in $X$.
When $d \in \{1,2,3\}$ these deserve a separate notation:
we write $\LX$ (resp. $\CX$, $\TX$) for the sets of divisor classes corresponding
respectively to lines, conics and twisted cubics in $X$.
There are respectively $27$, $27$ and $72$ of them, cfr. the appendix \ref{app-1}.

In the following proposition we will show that aCM line bundles on $X$ correspond to these divisor classes.
Although it is easy to classify aCM line bundles on any Del Pezzo surface
by geometric methods, we will outline an algebraic approach,
which is suitable for rank $2$ bundles as well.

\begin{prop} \label{acm-line-bundles}
Let $\LL$ be a normalized aCM line bundle on $X$.
Then the minimal graded free resolution of $\LL$ takes one of the
following forms:
\begin{small}
\begin{align}
  \label{res-OO}  & 0 \rr \OO(-3) \xr{F=\fff(\LL)} \OO \rr \LL \rr 0 & c_1(\LL) = 0\\
  \label{res-line}  & 0 \rr \OO(-2)^2 \xr{\fff({\LL})}
  \OO(-1) \oplus \OO \rr {\LL} \rr 0  & c_1(\LL) = L \in \LX \\
  \label{res-conic} & 0 \rr \OO(-1) \oplus \OO(-2) \xr{\fff({\LL})}
  \OO^2 \rr {\LL} \rr 0 & c_1(\LL)=C\in \CX \\
  \label{res-cubic} & 0 \rr \OO^3(-1) \xr{\fff({\LL})}
  \OO^3 \rr \LL \rr0 & c_1(\LL)=T \in \TX
\end{align}
\end{small}

Conversely, these divisor classes are associated to aCM line bundles.
In particular, there are $27$ (resp, $27$, $72$) ways
of writing $F$ as a determinant of the form 
(\ref{res-line}) (resp. (\ref{res-conic}), (\ref{res-cubic})).
\end{prop}

\begin{proof}
Recall the notation from Theorem \ref{res-general}.
Clearly $\rk(\gnrt(\LL))=1$ implies $\LL \simeq \OO_Y$, for $\ppp(\LL)_{|X}$
is surjective, hence an isomorphism. So, assume $\rk(\gnrt(\LL))\geq 2$.

By the minimality of the resolution, any degree zero term in the matrix $\fff(\LL)$ vanishes.
Thus, any summand contributing to the development of $\det(\fff(\LL))$
is either given by a product of three
linear forms of by a product of a quadratic form and a linear form.
Then the rank of $\gnrt(\LL)$ and $\syz(\LL)$ is either $3$ or $2$.
Furthermore, since $F$ is irreducible, any row and any column of the
the matrix $\fff(\LL)$ contains at least two nonvanishing entries.

Therefore if $\rk(\gnrt(\LL))=3$ all entries of $\fff(\LL)$ are linear
and the resolution takes the form (\ref{res-cubic}).
On the other hand if $\rk(\gnrt(\LL))=2$
the two summands in the development of $\det(\EE)$ are both a product
of a quadric and a linear form.
This gives cases (\ref{res-line}) and (\ref{res-conic}).
The remaining statements are clear.
\end{proof}

The following remark summarizes some of the classical combinatorics of these divisor classes.
Its proof is easy but rather tiresome: we leave to the contientious reader the task of verifying it,
making use of the tables in the Appendix \ref{app-2}.

\begin{rmk} \label{combinatorics}
Let $L,L_1,L_2 \in \LX$, $C,C_1,C_2 \in \CX$, $T,T_1,T_2 \in \TX$.
We have the well-defined surjective maps:
\begin{tiny}
\begin{align}
& \label{C+L=T} 
\{(T,C)  \,|\, T \cdot C = 1 \}  \xr{16:1} {\sf L}(X) && (T,C)  \mapsto T-C \\
& \label{L+C=T} 
\{(T,L)  \,|\, T \cdot L = 0 \}  \xr{16:1} {\sf C}(X) && (T,L)  \mapsto T-L \\
& \nonumber \{(C,L)  \,|\, C \cdot L = 0 \}  \xr{10:1} {\sf L}(X) && (C,L)  \mapsto C-L \\
& \nonumber \{(L_1,L_2)  \,|\, L_1 \cdot L_2 = 1 \} \xr{10:1} {\sf C}(X) && (L_1,L_2) \mapsto L_1+L_2 \\
& \nonumber \{(L_1,L_2)  \,|\, L_1 \cdot L_2 = 0 \} \xr{6:1} {\sf T}(X) && (L_1,L_2) \mapsto H-L_1+L_2 \\
& \label{twisted-lines} \{(T,L) \,|\, T \cdot L = 1 \} \xr{5:1} \{\{L_1,L_2\} \,|\, L_1 \cdot L_2 = 0 \} && (T,L) \mapsto T-L \\
& \label{three-lines} \{\{L_1,L_2,L_3\} \,|\, L_i \cdot L_j = 0 \} \xr{1:1} \{\{T_1,T_2\} | T_1 \cdot T_2 = 2 \} && L_1+L_2+L_3 \mapsto H+L_1+L_2+L_3
	 \\
& \nonumber \tau: {\sf T}(X)  \longleftrightarrow {\sf T}(X) && T  \mapsto 2\,H - T \\
& \nonumber \rho: {\sf L}(X)  \longleftrightarrow {\sf C}(X) && L  \mapsto H - L  
\end{align}
\end{tiny}
where the number over the arrow denotes the cardinality of the fibre.
\end{rmk}

\begin{rmk}
Assume $C+L=H$, i.e. $C=\rho(L)$.
Then the transpose of $\fff(\OO_X(L))$ 
(resp. of $\fff(\OO_X(C))$) gives a minimal resolution over $\p^3$ of $\OO_X(C+2\,H)$
(resp. of $\OO_X(L+2\,H)$).
Moreover, once restricted to $X$, 
we get the infinite $2$-periodic exact sequence:
\begin{footnotesize}
$$
\cdots 
\xr{\fff(\OO(L-3\,H))} \,\,
\begin{array}{c}
\OO(-3) \\
\oplus \\
\OO(-4)
\end{array}
  \xr{\fff(\OO(C-2\,H))} \,\, \OO(-2)^2 
 \xr{\fff(\OO(L))} \,\, 
\begin{array}{c}
\OO \\
\oplus \\
\OO(-1)
\end{array}
 \rr \,\, \OO(L)  \rr \,\, 0
$$
\end{footnotesize}

Similarly, $\fff({\OO_X(T+H)})^{\top}$ gives a resolution
over $\p^3$ of $\OO_X(2\,H - \tau(T))$.
\end{rmk}

Here we collect some elementary remarks about line bundles over $X$.
According to \cite[Theorem V.4.11]{hartshorne:ag}, a divisor class $D$ on $X$ is very ample iff
it is ample, iff it satisfies $D^2>0$ and $D \cdot L >0$, for $L\in \LX$.
By \cite[Exercise V.4.8]{hartshorne:ag} $D$ contains an integral curve iff it contains a smooth irreducible one,
iff it is a line, a conic, or $D$ satisfies $D^2>0$ and $D\cdot L\geq 0$ for all $L\in \LX$.

\begin{lem} \label{trivial}
Let $\LL\neq \OO_X$ be a line bundle on $X$ with $\hh^0(\LL) >0$, and let $\CC$ be an element in $|\LL|$.
Then the following hold.
\begin{enumerate}[i)]
\item \label{trivial-1}  If $\CC$ is reduced and connected, then $\HH^1(\LL(t-1))=0$ for $t\geq 0$.
\item \label{trivial-2}  If $\HH^1(\LL^*)=0$ then we have $\Ext^1(\OO_\CC,\LL^*)\simeq \kk$,
 and the unique extension class corresponds to the exact sequence defining $\CC \subset X$.
\item \label{trivial-0-ter} If $\CC$ is rational irreducible, then there exists $L \in {\sf L}(X)$ with $\LL \cdot L=0$.
\end{enumerate}
\end{lem}

\begin{proof}
Taking a section $s\in \HH^0(\LL)$ corresponding to $\CC$, we can write the two
equivalent exact sequences:
  \begin{align} 
    & \label{section-LL-1} 0 \rr \OO_X \xr{s} \LL \rr \OO_{\CC}(\CC
    \cdot \LL) \rr 0 \\
    & \label{section-LL-2} 0 \rr \LL^* \xr{s} \OO_X
    \rr \OO_{\CC} \rr 0
  \end{align}

Given a reduced connected curve $\CC$, we have $\HH^1(\OO_X(\CC)^*)=0$.
Thus, $\HH^1(\LL^*(-t))^*=\HH^1(\LL(t-1))=0$,
for $t \geq 0$, indeed a general curve in $|\LL(t)|$ is also reduced and connected.
So we have (\ref{trivial-1}).

For (\ref{trivial-2}), just apply $\Hom(-,\LL^*)$ to the sequence (\ref{section-LL-2}) defining $\CC$,
and observe that the image of the identity in $\End(\LL^*) \simeq \kk$ is the extension corresponding to (\ref{section-LL-2}) itself.

To check (\ref{trivial-0-ter}), notice that if $\CC$ is rational then we have $\LL^2 = \deg(\LL) - 2= \hh^0(\LL)-2$.
The statement is clear if $\deg(\LL) \leq 2$.
If $\deg(\LL) \geq 3$, then $\LL^2 \geq 1$.
Therefore, assuming $\LL \cdot L >0$ for $L\in \LX$, we deduce that $\LL$ is a very ample
line bundle.
Observe that $\LL$ would turn $X$ into a nondegenerate
surface of degree $m$ in $\p^{m-1}$, with rational hyperplane sections.
This is impossible (cfr. \cite[Exercise IV.4]{beauville:complex-algebraic-surfaces}).
\end{proof}

We will need to analyze linear systems containing rational curves.
The set of irreducible components of a reduced curve $\CC$ can be depicted via its intersection graph,
with one vertex for each component, and $n$ edges between two vertices if the two components meet at $n$ points.
This graph is connected iff $\CC$ is connected.

Given an effective divisor class $D$ on $X$, write $D \lequiv M+B$,
where $B$ is the {\em fixed component} of $|\OO_X(D)|$ (i.e. $B \subset \CC$, for all $\CC$ in $|\OO_X(D)|$)
and $M$ is the
{\em moving part}, (i.e. $|\OO_X(M)|$ has no fixed components).

\begin{lem} \label{elementary} Let $D$, $B$, $M$ be as above. Then:
\begin{enumerate}[i)]
\item \label{base-lines} there are disjoint lines $L_j \subset X$, {j=1,\ldots,b}, and integers $m_j \geq 0$ with
$B \lequiv  m_1 \, L_1 + \cdots + m_b \, L_b$;
\item \label{moving-and-base} we have $B \cdot M = 0$, i.e. a curve in $|\OO_X(M)|$ does not meet $B$;
\item \label{general-smooth} a general curve in $|\OO_X(M)|$ is smooth;
\item \label{general-irreducible}
a general curve in $|\OO_X(M)|$ is irreducible, unless $M \lequiv m\,C$, for some conic
$C \subset X$, and $m\geq 2$.
\end{enumerate}
\end{lem}

\begin{proof}
The proofs are based on the remarks above.
For (\ref{base-lines}), we write
$B \lequiv \sum_{j=1,\ldots,b} m_j B_j$, for some integral smooth curves $B_j \subset X$ and $m_j \geq 1$.
If $B_j^2>0$, or if $B_j$ is a conic, then $B_j$ cannot be a fixed component.
So $L_j$ must be a line. For $i \neq j$, $L_i$ does not meet $L_j$, for otherwise $L_i \cup L_j$
is linearly equivalent to a smooth conic, which moves in a pencil.

To show (\ref{moving-and-base}), 
choose $\CC$ to be a component of a curve in $|\OO_X(M)|$ and
let $L_j$ be a line in the support of $B$. Since $\HH^1(\OO_X(\CC))=0$, and since
the linear systems $|\OO_X(\CC)|$ and $|\OO_X(\CC + L_j)|$ have the same dimension
(indeed $L_j$ is a fixed component), we conclude $\CC \cdot L_j = 0$,
and so $M \cdot B = 0$.

The same argument implies $\CC \cdot L \geq 0$, for any line $L$. 
Now, given two components $\CC_1$, $\CC_2$ of a curve in $|\OO_X(M)|$,
if $\CC_1 \cdot \CC_2 > 0$, we have $(\CC_1 + \CC_2)^2 > 0$, so $\CC_1 + \CC_2$ is linearly equivalent
to a smooth integral curve. This proves (\ref{general-smooth}). For (\ref{general-irreducible}),
if $(\CC_1 \cdot \CC_2)^2 = 0$ we have $\CC_1^2 = \CC_2^2 = 0$, so the class of $\CC_i$ lies in $\CX$.
But then $\CC_1 \lequiv \CC_2$.
\end{proof}

\begin{lem} \label{other-splitting}
Let $\LL$ be a nontrivial line bundle on $X$ with $\hh^0(\LL) >0$ and
$\hh^1(\LL^*)\leq 1$. Let $\CC$ be a general curve in $|\LL|$.

\begin{enumerate}[i)]
\item \label{smooth} If $\hh^1(\LL^*)=0$, then $\CC$ is smooth irreducible.
\item \label{two-curves} If $\hh^1(\LL^*)=1$, then $\CC$ is the union of two disjoint smooth irreducible curves
$\CC_i \in |\LL_i|$, $i=1,2$ with $\LL \simeq \LL_1 \otimes \LL_2$. In this case there is (up to scalars) a unique
nonsplitting extension:
\begin{equation} \label{splitting}
0 \rr \LL^* \rr \LL_1^* \oplus \LL_2^* \rr \OO_X \rr 0
\end{equation}

In particular any extension corresponding to $\Ext^1(\OO_X,\LL)$ is a decomposable bundle.
\end{enumerate}
\end{lem}

\begin{proof}
Taking a nonzero section $s \in \HH^0(\LL)$ we obtain the exact
sequence (\ref{section-LL-2}) associated to the curve $\CC \in |\LL|$.
Write $c_1(\LL)$ as $B+M$, according to Lemma \ref{elementary}.
Then the number of connected components of a curve $\DD$ in $|\OO_X(M)|$ is $\hh^1(\OO_X(-M))+1$.
Notice that if $\hh^1(\OO_X(-M))=1$, $\DD$ must be the union of two linearly equivalent conics.
So, assume $B$ is nonempty. We write the exact sequence:
\begin{equation}
  0 \rr \LL^* \rr \OO_X(M) \rr \OO_{B} \rr 0
\end{equation}

Now, one checks easily that $\hh^0(\OO_{m \, L})=1/2(m^2+m)$ for $m\geq 1$.
So, if $M$ is empty, $\HH^1(\LL^*)=0$ implies that $B$ is a simple line,
while $B$ is the union of two skew lines if $\hh^1(\LL^*)=1$.

On the other hand, if $M$ is nonempty, we get $\hh^0(\OO_B) \geq \hh^1(\LL^*)$. 
Thus we get $\hh^1(\LL^*)\leq 1$, and (\ref{smooth}) is proved.
If $\hh^1(\LL^*)=1$, we conclude that $B$ is a simple line and $\HH^1(\OO_X(-M))=0$,
so $M$ is irreducible.

So in case (\ref{two-curves}) we have $\CC = \CC_1 \cup \CC_2$,
with $\CC_i$ given by the section $s_i$ of the a line bundle $\LL_i$, $i=1,2$.
Since $\CC_1 \cdot \CC_2=0$, the exact sequence (\ref{splitting})
is the Koszul complex of the section $(s_1,s_2) \in \HH^0(\LL_1 \oplus \LL_2)$.
But by $\HH^1(\LL^*) \simeq \Ext^1(\OO_X,\LL^*)=\kk$, one sees that such extension is unique.
\end{proof}

\begin{lem}
Take $R_d \in \RdX$. Then:
\begin{enumerate}[i)] \label{codimension}
\item \label{codimension-3} the set of nonreduced curves in $|\OO_X(R_d)|$ has codimension at least three;
\item \label{graph} a reduced curve in $|\OO_X(R_d)|$ is a simply connected graph of smooth rational curves;
\item \label{codimension-1} the set of reducible curves in $|\OO_X(R_d)|$ has codimension one.
\end{enumerate}
\end{lem}

\begin{proof}
A nonreduced component of a curve $\CC \in |\OO_X(R_d)|$ belongs to $|\OO_X(2\,R_e)|$, with
$R_e \in {\sf R}_e(X)$, for some $e \leq 2\,d$.
Since $\CC$ is connected, we have $(R_d-2\,R_e) \cdot R_e \geq 2$.
One proves easily that $\HH^1(\OO_X(R_d-2\,R_2))=\HH^2(\OO_X(R_d-2\,R_2))=0$, so
$\hh^0(\OO_X(R_d-2\,R_2))=d+e-4-R_d \cdot R_e \leq d-e-2$.
So nonreduced curves belong to subsets of the form $\p^{e-1} \times \p^{d-e-3}$, which have
codimension three in $\p(\OO_X(R_d))=\p^{d-1}$. This proves (\ref{codimension-3}).

For (\ref{graph}), notice that any component of a reduced curve in $|\OO_X(R_d)|$ is linearly equivalent
to a smooth rational curve, so it is itself smooth. The graph is simply connected since the arithmetic genus is zero.
Finally, by Lemma \ref{trivial}, (\ref{trivial-0-ter}), there is an $L\in \LX$ with
$R_{d-1}:=R_d-L \in {\sf R}_{d-1}(X)$, and $\p(\OO_X(R_{d-1}))$ is a codimension
one subset of reducible curves in $\p(\OO_X(R_d))$. This proves (\ref{codimension-1}).
\end{proof}

In the following Lemma we classify line
bundles of degree up to $3$ whose first cohomology group is trivial.

\begin{lem} \label{classify-deg-5}
Let $\LL$ be a line bundle on $X$ with $1 \leq \deg(\LL) \leq 3$, and
$\hh^1(\LL)=0$. Suppose $\hh^0(\LL)>0$ and let $\CC$ be a curve in
$|\LL|$. Then we have the following cases:

\begin{scriptsize}
\begin{center}
$$\newcounter{classify-LL}
\newcounter{sub-classify-LL}[classify-LL]
\begin{array}{c|c|c|c|c}
    \mbox{Ref.} & \deg(\LL) & c_1(\LL) & \hh^0(\LL) & \gen(\CC)  \\
    \hline 
    \hline
    \refstepcounter{classify-LL} \label{LL-1}
    \mbox{(L\arabic{classify-LL})} & 1 & L & 1 & 0  \\
    \hline
    \refstepcounter{classify-LL} \label{LL-2-1}
    \mbox{(L\arabic{classify-LL})} & 2 & C & 2 & 0  \\
    \refstepcounter{classify-LL} \label{LL-2-2}
    \mbox{(L\arabic{classify-LL})} & 2 & L_1+L_2 & 1 & -1  \\
    \hline
    \refstepcounter{classify-LL} \label{LL-3-1}
    \mbox{(L\arabic{classify-LL})} & 3 & H & 4 & 1  \\
    \refstepcounter{classify-LL} \label{LL-3-2}
    \mbox{(L\arabic{classify-LL})} & 3 & T & 3 & 0  \\
    \refstepcounter{classify-LL} \label{LL-3-3}
    \mbox{(L\arabic{classify-LL})} & 3 & C+L & 2 & -1  \\
    \refstepcounter{classify-LL} \label{LL-3-4}
    \mbox{(L\arabic{classify-LL})} & 3 & L_1+L_2+L_3 & 1 & -2 
\end{array}$$
\end{center}
\end{scriptsize}
where $L,L_i \in \LX$, $C\in \CX$, $T\in \TX$,  $L_i \cdot L_j= L \cdot C = 0$ for $i \neq j$.

Moreover, if $\hh^0(\LL)=0$, and $\hh^0(\LL(1))\neq 0$, then a general
curve $\DD$ in $|\LL(1)|$ is smooth, connected and rational of degree $3+\deg(\LL)$.
In this case we have:
\[
\hh^0(\LL(1))= \chi(\LL(1)) = 3 + \deg(\LL) \qquad \LL^2 = d+1
\]
\end{lem}

\begin{proof}

Recall the exact
sequence (\ref{section-LL-1}).
The case (L\ref{LL-1}) is obvious.
If $\deg(\LL)\in \{2,3\}$, our statement is equivalent to 
the claim that $\CC$ is reduced.
Given an integer $m$, consider the exact sequence:
  \begin{equation}
    \label{2L} 0 \rr \OO_X((m-1)\,L) \xr{s} \OO_X(m\,L) \rr \OO_{L}(-m) \rr 0
  \end{equation}

Setting $m=2$ (resp. $m=3$) in (\ref{2L}), we see that
$\hh^1(\OO_X(2\,L)) = 1$ (resp. that $\hh^1(\OO_X(3\,L)) = 3$.
So the fixed components of $\CC$ are reduced, and we are done with the first claim.

Now suppose $\hh^0(\LL)=0$, $\hh^0(\LL(1))\neq 0$, and let $\DD$ be a
general curve of the linear system $|\LL(1)|$.
By Serre duality $\HH^1(\OO_X(\DD)^*)=0$, so by Lemma \ref{other-splitting} $\DD$ is smooth irreducible.
Now $\hh^2(\LL)=0$ implies that $\DD$ is rational. The last formulas follow easily.
\end{proof}

Next we classify all nontrivial extensions of two aCM line bundles $\MM$ and $\NN$ on $X$.
For the next lemma, we set $M=c_1(\MM)$, $N = c_1(\NN)$, $\Sigma = c_1(\MM \ts \NN)$,
$\Delta = c_1(\MM \ts \NN^*)$,
$\hh^1$ will indicate the dimension of $\HH^1(\MM \ts \NN^*(t))$, $\delta = (\Delta+H)^2$,
$\sigma = (\Sigma-H)^2$. Here $R_d$ will stand for an element of $\RdX$.
Notice that the integers $\delta,\sigma$ determine the intersection of the divisors appearing in the
expression of $\Delta + H$ and $\Sigma-H$.

\begin{lem} \label{all-extensions}
Set notations as above. The
group $\HH^1(\MM \ts \NN^*(t))$ vanishes except in the
cases comprehended by the following table.

\begin{tiny}
\begin{align}
\label{T-T}
\fbox{$M\in\TX$} \quad &
\left\{
\begin{tabular}{c|c|c|c|c|c|c|c}
$N$ & $M \cdot N$ & $t$ & $\hh^1$ & $\Delta+H$ & $\delta$ & $\Sigma-H$ & $\sigma$ \\
\hline
\hline
$T$ & 5 & $(-1,0)$ & $(3,3)$ & $2\,M-H$ & $-5$ & $H$ & $3$ \\
$T$ & 4 & $(-1,0)$ & $(2,2)$ & $L_1+L_2+L_3$ & $-3$ & $T_1$ & $1$ \\
$T$ & 3 & $(-1,0)$ & $(1,1)$ & $C_1+L_1$ & $-1$ & $C_2+L_2$ & $-1$\\
\hline
$C_1$ & 3 & $(-1,0)$ 	& $(2,1)$ & $T+L$ & $0$ &  $C_2$ & $0$\\
$C$ & 2 & $-1$ 			& $1$ 		  & $R_4$ & $2$ & $L_1+L_2$ & $-2$ \\
\hline
$L_1$ & 2 & $(-2,-1)$ & $(1,2)$ & $R_5$ & $3$ & $L_2$ & $-1$\\
$L$ & 1 & $-1$ & $1$ & $L_1+L_2+H$ & $5$ & $H-L_3-L_4$ & $-3$\\
\end{tabular}
\right.\\
\label{C-C}
\fbox{$M \in \CX$} \quad &
\left\{
\begin{tabular}{c|c|c|c|c|c|c|c}
$N$ & $M \cdot N$ & $t$ & $\hh^1$ & $\Delta+H$ & $\delta$ & $\Sigma - H$ & $\sigma$ \\
\hline
\hline
$T$ & 3 & $(-1,0)$ & $(1,2)$ & $R_5-H$ & $-4$ & $C$ & $0$ \\
$T$ & 2 & $0$ & $1$ & $L_1+L_2$ & $-2$ & $L_3+L_4$ & $-2$ \\
\hline
$C_1$ & 2 & $(-1,0)$ & $(1,1)$ & $C_2 + \rho(C_1)$ &$0$ & $L$ & $-1$ \\
\hline
$L$ & 2 & $(-2,-1,0)$ & $(1,2,1)$ & $R_4-H$ & $-4$ & $0$ & $0$ \\
$L$ & 1 & $-1$ & $1$ & $L_1+L_2$ & $-2$ & $T-H$ & $-2$ \\
\end{tabular} \right. \\
\label{L-L}
\fbox{$M\in \LX$} \quad &
\left\{
\begin{tabular}{c|c|c|c|c|c|c|c}
$N$ & $M \cdot N$ & $t$ & $\hh^1$ & $\Delta+H$ & $\delta $ & $\Sigma - H$ & $\sigma$ \\
\hline
\hline
$T_1$ & 2 & $(0,1)$ & $(2,1)$ & $T_2+L-H$ & $-5$ &  $L$ & $-1$ \\
$T$ & 1 & $0$ & $1$ & $R_4-H$ & $-3$ & $H-L_1-L_2$ & $-3$ \\
\hline
$C$ & 2 &  $(-1,0,1)$ & $(1,2,1)$ & $L+R_4-H$ &$-4$& $0$ & $0$ \\
$C$ & 1 &  $0$ & $1$ & $L_1+L_2$ & $-2$& $T-H$ & $-2$ \\
\hline
$L_1$ & 1 & $(-1,0)$ & $(1,1)$ &$L_2+C$ & $-1$ & $-L_3$ & $-1$
\end{tabular}
\right.
\end{align}
\end{tiny}
\end{lem}

\begin{proof}
We classify the $\Delta$'s which are not aCM divisors.
If $M=N$, there is nothing to prove.
By the symmetry, we can assume $\deg(M) \leq \deg(N)$.

Consider the divisor $\Delta + H$. We have $1 \leq \deg(\Delta + H) \leq 3$.
By Lemma \ref{trivial}, (\ref{trivial-1}), $\HH^1(\OO_X(\Delta + t\,H))$ vanishes
for all $t\geq 1$ if $|\OO_X(\Delta + 2\,H)|$ contains reduced connected curves.
One proves immediately that this holds in all cases except $N \in \TX$, $M \in \LX$, $M \cdot N = 2$,
and $M = \rho(N) \in \LX$. It is easy to study these two cases separately.
Namely they give respectively $\Delta + 2\,H = M+\tau(N)$, and $\Delta + 2\,H = M+R_4$,
which are easily handled. Now one checks the formulas:
\begin{small}
\begin{align}
\label{chi-sigma} & \hh^0(\OO_X(\Sigma-H))=\chi(\OO_X(\Sigma-H)) = M \cdot N - 1 \\
\label{chi-delta} & \hh^0(\OO_X(\Delta+H))=\chi(\OO_X(\Delta+H)) = 2\,\deg(M)-\deg(N) + 2 - M \cdot N
\end{align}
\end{small}

Apply Lemma \ref{classify-deg-5} to $\OO_X(\Delta + H)$.
For positive values of (\ref{chi-delta}), we get a classification of $\Delta$.
If $\chi(\OO_X(\Delta+H))=0$, 
the description of the linear system $|\OO_X(\Delta + 2\,H)|$
is summarized by the table (where $\chi$ denotes $\chi(\OO_X(\Delta+H))$):
\begin{center}
\begin{tiny}
\begin{tabular}{c|c|c|c|c|c|c}
$\deg(M)$ & $\deg(N)$ & $M\cdot N$ & $\deg(\Delta + 2\,H)$ & $\chi$ & $\Delta+2\,H$ & $(\Delta+2\,H)^2$ \\
\hline
\hline
3 & 3 & 5 & 6 & 6 & $2\, \tau(N) = R_6$ & 4\\ 
\hline
2 & 3 & 3 & 5 & 5 & $ M + \tau(N) = R_5$ & 3\\ 
\hline
1 & 3 & 1 & 4 & 4 & $M + \tau(N) = R_4$ & 2\\ 
1 & 3 & 2 & 4 & 3 & $M + \tau(N) = T + L$ & 0\\ 
\hline
1 & 2 & 2 & 5 & 4 & $M+H+\rho(N)=M+R_4$ & 1
\end{tabular}
\end{tiny}
\end{center}

So we have classified the $\Delta$'s. We see that $\HH^1(\OO_X(\Delta - t\,H))=0$ for $t\geq 2$, since in all cases
$|\OO_X(2\,H-\Delta)|$ contains reduced connected curves. Summing up we have $t\in \{-1,0\}$
(except in the two cases treated separately).
We leave now to the reader the exercise of computing the value of $\hh^1(\OO_X(\Delta + t\,H))$.

To finish the proof, it remains to compute $\Sigma$. We observe:
\begin{footnotesize}
\[
\Delta + t\,H = \left\{						\begin{array}{ll}
																		M+\rho(N) + (t-1)\,H & \mbox{if $N \in \LX \cup \CX$} \\
																		M+\tau(N) + (t-2)\,H & \mbox{if $N \in \TX $}
																	 	\end{array}
																	 	\right.
\]
\end{footnotesize}

According to the above alternative, set $\Delta'=M-N'$ with $N' := \tau(N)$ or $N' := \rho(N)$.
In the former case, we get $M \cdot N' = 2\,\deg(M)-M\cdot N$, 
and $\Sigma = M + 2\,H-N'=\Delta' + 2\,H$. In the latter case, we have
$M \cdot N' = \deg(M)-M\cdot N$ and $\Sigma = M + H-N'=\Delta' + H$.
But we have already classified $\Delta'$. 
\end{proof}

\section{Classification of resolutions} \label{clas-res}

Let $\EE$ will be a rank $2$ indecomposable aCM bundle on $X$.
We classify the degree of the generators $\gnrt(\EE)$
and syzygies $\syz(\EE)$ appearing in the
minimal graded free resolution of $\EE$, extended to zero to $\p^3$,
according to Theorem \ref{res-general}.

In the following theorem,
the column {\em dual} describes the minimal graded free resolution of $\EE^*$.
The column {\em kernel} provides a resolution (possibly nonminimal) of the aCM
vector bundle $\ker(\ppp(\EE)_{|X})$, in the case that it also has rank $2$
(i.e. in case $\rk(\gnrt(\EE))=4$).
In these two columns, the number in parenthesis points out the twist
in which the resolution of $\EE^*$ or $\ppp(\fff(\EE)_{|X})$ occurs.
In the {\em Hilbert} column we write the Hilbert polynomial of $\EE$.

\begin{thm}  \label{full-classification}
Let $X$ and $\EE$ be as above.
Then the minimal graded free resolution of $\EE$ takes one of the following forms:

\begin{minipage}[h]{11cm}

\begin{enumerate}[(A)]

\item $\gnrt(\EE) = \OO^6$ and $\syz(\EE)=\OO(-1)^6$; \label{6}
\item $\gnrt(\EE) = \OO^5$ and $\syz(\EE)=\OO(-1)^4 \oplus \OO(-2)$; \label{5-2}
\item $\gnrt(\EE) = \OO \oplus \OO(-1)^4$ and $\syz(\EE)=\OO(-2)^5$; \label{5-4}
\item $\gnrt(\EE) = \OO^4$ and $\syz(\EE)=\OO(-1)^2 \oplus \OO(-2)^2$; \label{4-1}
\item $\gnrt(\EE) = \OO^3 \oplus \OO(-1)$ and $\syz(\EE)=\OO(-1)
  \oplus \OO(-2)^3$; \label{4-3}
\item $\gnrt(\EE) = \OO^2 \oplus \OO(-1)^2$ and $\syz(\EE)=\OO(-2)^4$;
  \label{4-2}
\item $\gnrt(\EE) = \OO \oplus \OO(-1)^3$ and $\syz(\EE)=\OO(-2)^3
  \oplus \OO(-3)$;\label{4-4} 
\item $\gnrt(\EE) = \OO^3$ and $\syz(\EE)=\OO(-2)^3$. \label{3}
\end{enumerate}

\end{minipage}

\vspace{0.3cm}
\noindent Moreover we can summarize the following information:
\begin{footnotesize}
\begin{equation} \label{summary} 
\begin{tabular}{c|c|c|c|c|c}
    Ref. & $\rk(\gnrt(\EE))$ & $\deg(c_1(\EE))$ & Hilbert & dual & kernel\\
    \hline
    \hline
    (\ref{6}) & 6 & 6  & $3\,t^2 +9\,t +6$ & \ref{6}(-2) &  \\
    \hline
    (\ref{5-2}) & 5 & 5 & ${3}\,t^2 +{8}\,t + 5$  & (\ref{5-4})(-1)  &  \\
    (\ref{5-4}) & 5 & 1 & ${3}\,t^2 +4\,t + 1$ & (\ref{5-2})(-1) &  \\
    \hline
    (\ref{4-1}) & 4 & 4 & ${3}\,t^2 +7\,t + 4$ & (\ref{4-2})(-1) & (\ref{4-2})(-1) \\
    (\ref{4-3}) & 4 & 3 & ${3}\,t^2 +6\,t + 3$ & (\ref{4-3})(-1) & (\ref{4-4})(-1) \\
    (\ref{4-2}) & 4 & 2 & ${3}\,t^2 +5\,t + 2$ & (\ref{4-1})(-1) & (\ref{4-1})(-2) \\
    (\ref{4-4}) & 4 & 0 & ${3}\,t^2 +3\,t + 1$ & (\ref{4-4}) & (\ref{4-3})(-2) \\
    \hline
    (\ref{3}) & 3 & 3 & ${3}\,t^2 +6\,t + 3$ & (\ref{3})(-1) & 
\end{tabular}
\end{equation}
\end{footnotesize}
\end{thm}

\begin{proof}[Beginning of the proof of \ref{full-classification}]
The proof is similar to that of Proposition \ref{acm-line-bundles},
though more involved. In view of Theorem \ref{res-general} we consider the matrix
$\fff(\EE)$ in the 
minimal graded free resolution of ${\EE}$,
satisfying $\det(\fff(\EE))=F^2$.
Since we assume that the resolution is minimal,
any entry of degree zero in the matrix $\fff(\EE)$ vanishes.
Therefore, we have $\rk(\gnrt(\EE)) \leq 6$.

Clearly, we have $\rk(\gnrt(\EE))\geq 2$, and equality holds if and
only if $\EE$ is isomorphic to $\OO_X \oplus \OO_X(-m)$, for some
$m\geq 0$.
Indeed if $\rk(\gnrt(\EE)) = 2$, the map $\ppp(\EE)_{|X}$ is a surjective morphism of
vector bundles of the same rank, hence an isomorphism.

We split the proof into cases, according to $\rk(\gnrt(\EE))=3$, $4$, $5$, $6$.
\end{proof}

\begin{lem}
Let $\EE$ and $X$ be as in Theorem \ref{full-classification}, and
suppose $\rk(\gnrt(\EE))=3$.
Then $\fff(\EE)$ is a matrix of quadratic forms, i.e. $\EE$ is of type (\ref{3}).
Moreover, there is $T\in \TX$ such that:
 \begin{equation} \label{wedge-2} \fff(\EE) = \wedge^2 \fff(\OO_X(T))
 \end{equation}
\end{lem}

\begin{proof}
Set $\LL := \ker(\ppp(\EE)_{|X})$.
By Remark \ref{periodicity}, $\LL$ is an aCM line bundle, and there exists a matrix $g$
defined on $\p^3$ such that $\coker(g) \simeq \LL$, $g_{|X} \circ
\fff(\EE)_{|X}=0$ and $\fff(\EE)_{|X} \circ g_{|X}=0$.

Recall now Proposition \ref{acm-line-bundles}.
In case $\rk(\gnrt(\LL))=3$, $\fff(\LL)$ is a matrix
of linear forms, and $\LL \simeq \OO_X(T)$, for some $T\in \TX$. 
Since $\det(g) = F$, by $g \circ \wedge^2(g) = F \, \one_{3}$ we obtain (\ref{wedge-2}).
Thus $\EE$ is of type (\ref{3}).

If $\LL$ is a twist of $\OO_X(L)$ or $\OO_X(C)$ (i.e. if
$\rk(\gnrt(\LL))=2$), the matrix $g$ can be reduced, under the action
by conjugation of the group $\GL(3,\kk)$, to:
$$g = \left(\begin{array}{c|c}
    0 & \lambda \\
    \hline
    g' & 0
    \end{array} \right) \quad \mbox{with $0 \neq \lambda \in \kk$, and where $g'$
    is a $2 \times 2$ matrix.}$$

Since $g_{|X} \circ \fff(\EE)_{|X}=0$ and $\fff(\EE)_{|X} \circ
g_{|X}=0$, $\fff(\EE)$ can be reduced to: 
$$\fff(\EE) = \left(\begin{array}{c|c}
    0 & F \\
    \hline
    f' & 0
    \end{array} \right) \quad \mbox{where $f'$ is a $2 \times 2$ matrix.}$$

This implies that $\EE$ is decomposes as $\OO_X(m) \oplus \coker(f')$,
for some $m$.
On the other hand, if $\rk(\gnrt(\LL))=1$ (i.e. if $\LL \simeq
\OO_X(m)$, for some $m$), the matrix $g$ can be reduced to:
$$g = \left(\begin{array}{c|c}
    0 & F \\
    \hline
    g' & 0
    \end{array} \right) \quad \mbox{where $g'$
    is a $2 \times 2$ invertible matrix.}$$

Also in this case, the bundle $\EE$ is decomposable.
\end{proof}

Now we assume $\rk(\gnrt(\EE))\geq 4$. 
Any nonzero summand contributing to the development of $\det(\EE)$
is given by a product of one of the following types:
\begin{enumerate}[a)]
\item \label{rank-6} Six linear entries;
\item \label{rank-5} Four linear entries and one quadratic entry;
\item \label{rank-4-I} Three linear entries and one cubic entry;
\item \label{rank-4-II} Two linear entries and two quadratic entries.
\end{enumerate}

Clearly the rank of $\gnrt(\EE)$ is determined by the above alternatives.
We will analyze separately the cases $\rk(\gnrt(\EE))=4,5,6$ in the
following lemmas. We need the following claim, which is an analogue of a result of Bohnhorst-Spindler \cite{bohnhorst-spindler}.

\begin{claim}
The sequences of integers $a_1,\ldots,a_s$ and $b_1,\ldots,b_s$ satisfy the following relations:
\begin{equation} \label{bohnhorst-spindler}
\sum_{i=1}^s b_i - \sum_{i=1}^s a_i = 6 \qquad a_{\ell} \leq b_{\ell} - 1 \,\, \mbox{for each $\ell=1,\ldots,s$}
\end{equation}
\end{claim}

\begin{proof}
Let $\fff = \fff(\EE)$.
The first claim is obvious since $\deg(\det(\fff))=6$.
Now fix an integer $\ell$ with $1 \leq \ell \leq 6$ and
consider the maximal number $r$ such that $(\fff_{r,1},\ldots,\fff_{r,\ell}) \neq (0,\ldots,0)$.
Notice that $r\geq \ell$ for otherwise $\fff$ is not injective.
Then $\fff_{i,j}$ gives an injective map:
$$\bigoplus_{i=1}^{\ell} \OO_{\p^3}(a_i) \rr \bigoplus_{j=1}^{r} \OO_{\p^3}(b_j)$$

Then $\fff_{r,i} \neq 0$ for some $i\leq \ell$, so $a_{\ell} \leq a_i < b_r - 1 \leq b_{\ell} -1$, for each $\ell$.
\end{proof}

\begin{lem}[Six by Six] \label{six-by-six}
Let $\EE$ and $X$ be as in Theorem \ref{full-classification}, and
suppose $\rk(\gnrt(\EE))=6$.
Then $\fff(\EE)$ is a matrix of linear forms.
\end{lem}

\begin{proof}
Set $\fff = \fff(\EE)$.
Arrange a set of entries of type (\ref{rank-6}) along the main diagonal of $\fff$.
Let $\beta_i = b_i-a_{i+1}$, and notice by Claim \ref{bohnhorst-spindler} that
$\beta_i \geq 1$ for each $i$. Write $\beta$ for the sequence of $\beta_i$'s.
We have to prove that $\beta_i = 1$, for $1 \leq i \leq 5$, i.e. $\beta = (1,1,1,1,1)$.

If $\beta_1 \geq 2$, then $\fff_{1,1}$ divides $\det(\fff)$, a contradiction.
The same happens if $\beta_5 \geq 2$, so $\beta_1 = \beta_5 = 1$.
Analogously, if $\beta_2 \geq 2$, the minor $\wedge^2(\fff)_{1,1}$ divides
$\det(\fff)$. But this minor has degree $2$ and $F$ is irreducible: a contradiction.
Similarly for $\beta_4 \geq 2$, so $\beta_2 = \beta_4 = 1$.

Assuming $\beta_3\geq 2$, the minor $\wedge^3(\fff)_{1,1}$ divides $\det(\fff)$.
Let then ${\sf g}^1$ and ${\sf g}^2$ be the two $3 \times 3$ submatrices of linear forms sitting on the main diagonal of $\fff$. The determinant of both ${\sf g}^1$ and ${\sf g}^2$ is a multiple of $F$,
so that 
$\cok({\sf g}^i)$ is a twist of $\OO_X(T_i)$ for some $T_1, T_2 \in \TX$.
In fact we can construct a commutative exact diagram:

\begin{small}
\begin{equation} \label{force-extension}
\xymatrix{
{\OO_{\p^3}(-1)^3} \ar^-{{\sf g}^1}[r] \ar[d] & {\OO_{\p^3}^3} \ar[r] \ar[d] &
{\OO_X(T_1)} \ar@{-->}[d] \\
{\OO_{\p^3}(-1)^3 \oplus \OO_{\p^3}(-\beta_3)^3} \ar^-{\fff}[r] \ar[d] & {\OO_{\p^3}^3 \oplus \OO_{\p^3}(1-\beta_3)^3} \ar[r] \ar[d] &
{\EE} \ar@{-->}[d] \\
{\OO_{\p^3}(-\beta_3)^3} \ar^-{{\sf g}^2}[r]  & {\OO_{\p^3}(1-\beta_3)^3} \ar[r] & {\OO_X(T_2+(1-\beta_3)\,H)}  \\
}
\end{equation}
\end{small}
where the solid vertical
maps are the inclusions and projections corresponding to the block subdivision of $\fff(\EE)$, the dashed
maps are induced on $X$, and we omit zeroes
all around the diagram for brevity.
But by Lemma \ref{all-extensions} we have:
\begin{footnotesize}
\[
\Ext^1(\OO_X(T_2+(1-\beta_3)\, H),\OO_X(T_1)) \simeq
\HH^1(\OO_X(T_1 - T_2 + (\beta_3-1)\, H)) =0
\]
\end{footnotesize}
for $\beta_3\geq 2$, and for any pair $(T_1,T_2)$, so $\EE$ is decomposable.
\end{proof}

\begin{lem}[Five by Five] \label{five-by-five}
Let $\EE$ and $X$ be as above, and
suppose $\rk(\gnrt(\EE))=5$.
Then $\fff(\EE)$ takes one of the forms (\ref{5-2}), (\ref{5-4}),
of Theorem \ref{full-classification}.
\end{lem}

\begin{proof}
Arrange a set of entries of type (\ref{rank-5}), ordered by ascending degree,
on the main diagonal of $\fff(\EE)$, and
recall the notation from the proof of Lemma \ref{six-by-six}.
We would like to prove that $\beta$ must take value $(1,1,1,1)$ or $(1,1,1,2)$.

If $\beta_1 \geq 2$, looking at the first column of $\fff$ we conclude by the irreducibility of
$F$ that the only possibility is $\beta_1 = (2,1,1,1)$.
We will show that this does not occur at the end of the proof.

So let's assume $\beta_1=1$, $\beta_2 \geq 2$. Here irreducibility of $F$ implies $\beta = (1,2,1,1)$.
This case also will be excluded at the end of the proof.

Then we suppose $\beta_1=\beta_2=1$, $\beta_3 \geq 2$.
Notice that $\beta_4 \geq 3$ contradicts irreducibility of $F$, so $\beta_3 \in \{1,2\}$.
Let us first look at $\beta_3=1$.

Since $\beta = (1,1,2,1)$ is equivalent to $\beta = (2,1,1,1)$ after transposition, we can assume $\beta_3 \geq 3$.
Thus we can construct a commutative diagram
analogous to (\ref{force-extension}), which implies 
that $\EE$ fits in the exact sequence:
$$0 \rr \OO_X(T+(\beta_3-2)\,H) \rr \EE \rr \OO_X(L) \rr 0$$
for some $L \in \LX$, $T\in \TX$. But Lemma \ref{all-extensions} implies:
\begin{footnotesize}
\[
\Ext^1(\OO_X(L),\OO_X(T+(\beta_3-2)\,H))=\HH^1(\OO_X(T-L+(\beta_3-2)\,H))=0
\]
\end{footnotesize}
for any pair $(T,L)$ and $\beta_3\geq 3$, so $\EE$ splits.
Similarly, if $\beta_4=2$ we conclude that $\EE$ splits by:
  $$\HH^1(\OO_X(T-C+(\beta_3-1)\,H))=0 \qquad \mbox{for any $C,T$ and for $\beta_3\neq 0,1$}$$

It remains to show that $\beta$ cannot be $(1,2,1,1)$ or $(2,1,1,1)$.
Considering the first case, 
after a permutation of the basis, we have to show that the following configuration does not occur:
\begin{tiny}
\begin{equation} \label{5-impossible}
\left(
\begin{array}{cc|ccc}
0 & 0 & 1 & 1 & 1\\
0 & 0 & 1 & 1 & 1\\
\hline
 1 & 1 & 2 & 2 & 2 \\
 1 & 1 & 2 & 2 & 2 \\
 1 & 1 & 2 & 2 & 2
\end{array}
\right)
\end{equation}
\end{tiny}

Write ${\sf g}^{i}$ (resp. $_i{\sf g}$) for the $2\times 2$ submatrix of the upper-right
(resp. lower-left) linear block of $\fff$ obtained deleting the $i$-th column
(resp. row), for $i\in \{3,4,5\}$. By Remark \ref{periodic}, part (\ref{minors}), since
we have,
for all $i,j$:
$$\det({\sf g}^i) \, \det(_j{\sf g}) = 0 \qquad \mbox{over $X$}$$

But since $\deg(_i{\sf g})=2$, we deduce $\det(_i{\sf g})=0$ over $\p^3$,
for all $i$ (unless the same happens to ${\sf g}^i$, for all $i$).
Hence $\det(\fff)=0$,
a contradiction.
To complete the proof, it remains only to exclude $\beta = (2,1,1,1)$,
i.e. the configuration:

\begin{tiny}
\begin{equation} \label{5-impossible-2}
\left(
\begin{array}{ccc|cc}
0 & 0 & 0 & 1 & 1\\
\hline
1 & 1 & 1 & 2 & 2\\
1 & 1 & 1 & 2 & 2\\
1 & 1 & 1 & 2 & 2 \\
1 & 1 & 1 & 2 & 2 \\
\end{array}
\right)
\end{equation}
\end{tiny}

Let ${\sf g}$ be the $4 \times 3$ submatrix of $\fff$
containing linear entries.
After Remark \ref{periodic}, part (\ref{minors}), $\wedge^3({\sf g})$
must vanish on $X$.
Nevertheless $\wedge^3({\sf g})$ cannot vanish identically on $\p^3$ for otherwise
$\det(\fff)=0$.
So, we can choose a $3\times 3$ submatrix of ${\sf g}$ whose determinant is a nonzero multiple of $F$.
In turn, its cokernel is isomorphic to $\OO_X(T)$, for some twisted cubic $T \subset X$.
Thus, we obtain that $\FF:=\coker({\sf g}_{|X})$
is a locally free sheaf on $X$ which fits into:
 $$0 \rr \OO_X \rr \FF \rr \OO_X(T) \rr 0$$

So $\FF$ splits for $\OO_X(T)$ is aCM.
Then we can assume that one row in the matrix ${\sf g}$ is zero.
In turn,  we obtain an exact sequence:
$$0 \rr \OO_X(T) \rr \EE \rr \OO_X(L) \rr 0$$
for some line $L \subset X$. So $\EE$ splits by Lemma \ref{all-extensions}.

\end{proof}

\begin{lem}[Four by four] 
Let $\EE$ and $X$ be as above, and
suppose $\rk(\gnrt(\EE))=4$.
Then $\fff(\EE)$ takes one of the forms (\ref{4-1}), (\ref{4-2}),
(\ref{4-3}) or (\ref{4-4}) of Theorem \ref{full-classification}.
\end{lem}

\begin{proof}
Recall the notation from the proof of the previous lemmas.
We divide the proof into two cases, according to the assumption that:
\begin{enumerate}[I)]
\item \label{proof-4-I} the matrix $\fff$ contains at least a
  set of entries of type (\ref{rank-4-I}); 
\item \label{proof-4-II} all summands in the development of
  $\det(\fff)$ are of type (\ref{rank-4-II}). 
\end{enumerate}

{\bf Case (\ref{proof-4-I})}.
Arrange the set of type (\ref{rank-4-I}) by ascending degree along the main diagonal of
$\fff$. We would like to show that $\beta = (1,1,2)$, i.e. case (\ref{4-4}).

Reasoning like in the proof of Lemma \ref{five-by-five}, one sees that
$\beta_1 \geq 2$ gives rise to the cases $\beta = (3,1,1)$, $(2,1,2)$, $(2,2,1)$, $(2,1,1)$.
In these cases, we
consider the $3\times 3$ submatrices of $\fff$
containing entries of the first column.
Since $\wedge^3(\fff)$ must vanish on $X$, it is easy to show that either $\fff_{1,1} = 0$,
either the remaining $2\times 2$ vanishes on $\p^3$.
In the latter case one writes a diagram similar to (\ref{force-extension})
and uses Lemma \ref{all-extensions} (\eqref{C-C}), (\eqref{L-L})
to conclude that $\EE$ is decomposable.

Now if $\beta_1 = 1$, $\beta_2 \geq 2$ gives rise to the cases:
$\beta = (1,2,1)$, $(1,3,1)$, $(1,2,2)$. The same argument as above goes through here.

So let us assume $\beta_1 = \beta_2 = 1$. If $\beta_3 \geq 4$, a diagram analogous to 
(\ref{force-extension}) shows that $\EE$ decomposes as $\OO_X(T) \oplus \OO_X(m)$, for some $T\subset X$
and some integer $m$.
Notice that $\beta = (1,1,3)$, and $\beta = (1,1,1)$ are equivalent after transposition, hence it only remains to
exclude $\beta = (1,1,3)$, i.e. the configuration:
\begin{tiny}
\begin{equation} \label{4-impossible-2}
\left(
\begin{array}{ccc|c}
  1 & 1 & 1 & 3 \\
  1 & 1 & 1 & 3 \\
  1 & 1 & 1 & 3  \\
  1 & 1 & 1 & 3 
\end{array}
\right)
\end{equation}
\end{tiny}

In this case we call ${\sf g}$ the $4 \times 3$ block of linear
entries in (\ref{4-impossible-2}).
We get that $\wedge^3({\sf g})$ vanishes on $X$,
and we conclude by the same argument as for (\ref{5-impossible-2}).

{\bf Case (\ref{proof-4-II})}.
We would like to prove that $\beta$ equals
$(1,1,2)$ (case (\ref{4-2})), $(1,2,2)$ (case (\ref{4-1})),
or $(2,1,2)$ (case (\ref{4-3})).

If $\beta_1 \geq 2$, we must only exclude $\beta = (2,2,1)$, but this is the same as $(2,1,2)$ of the case
(\ref{proof-4-II}) above. 

Assuming $\beta_1 = 2$, we see that it remains to exclude the following values of $\beta$:
$(1,3,1)$, $(1,2,1)$, $(1,2,3)$, $(1,1,3)$, $(1,1,1)$.
But one checks immediately that all these cases have already been taken into account,
up to transposition and permutation of the basis.
\end{proof}

\begin{proof}[End of the proof of \ref{full-classification}]
We have proved in the above Lemmas that the resolution of $\EE$ 
takes one of the desired forms.
This gives at once the Hilbert polynomial of $\EE$, and thus
the value of $\deg(\det(\EE))$.

By duality we have the formulas:
\begin{align}
\label{grothendieck} & \EExt^1(\EE,\OO_{\p^3}) \simeq \EE^*(3) \\
\label{omega-1} & \hh^1(X,\EE \ts \Omega^1_{\p^3}(1)) = \hh^1(X,\EE^* \ts \Omega^2_{\p^3}(2)) \\
\label{omega-2} & \hh^0(X,\EE \ts \Omega^1_{\p^3}(1)) = \hh^2(X,\EE^* \ts \Omega^2_{\p^3}(2)) \\
\label{serre-duality} & \hh^0(X,\EE(t))=\hh^2(X,\EE^*(-t-1))
\end{align}

Notice that, in order to determine the minimal
resolution of $\EE$ it suffices to 
compute its Hilbert function or its Hilbert polynomial, except when this equals
$3\,t^2+6\,t+3$. However, in this case we can distinguish between
(\ref{4-3}) and (\ref{3}) by:
\begin{align}
  \label{distinguish-1} & \hh^0(\EE \ts \Omega_{\p^3}(1)) = \hh^1(\EE \ts
  \Omega_{\p^3}(1))=0 && \mbox{in case (\ref{3})} \\
  \label{distinguish-2} & \hh^0(\EE \ts \Omega_{\p^3}(1)) = \hh^1(\EE \ts
  \Omega_{\p^3}(1))=1 && \mbox{in case (\ref{4-3})}
\end{align}

Therefore, in order to compute the minimal graded free resolution of $\EE^*$
it suffices to compute the Hilbert polynomial by means of (\ref{serre-duality}),
except in the cases (\ref{4-3}) and (\ref{3}).
But in these cases, making use of (\ref{omega-1}) and (\ref{omega-2}), the desired
resolution for $\EE^*$ is deduced by (\ref{distinguish-1}) and (\ref{distinguish-2}).
Finally, a resolution of the syzygy bundle $\ker(\ppp(\EE)_{|X})$ is
given by Remark \ref{periodicity}.
\end{proof}

\section{Chern classes of rank 2 aCM bundles} \label{chern}

We will classify the Chern classes of an indecomposable aCM bundle rank 2 bundle $\EE$ on $X$ according to
its minimal graded free resolution, cfr. Theorem \ref{full-classification}.
The next Theorem summarizes the results of this section.

\begin{thm} \label{full-chern}
Let $\EE$ be as above. Then the Chern classes of $\EE$ behave according to the following table.
\begin{scriptsize}
\begin{equation} \label{all-chern-classes}
\begin{array}{c|c}
 & \mbox{\em Chern} \\
\hline
\hline
\begin{array}{c}
\mbox{Ref.} \\
\hline
\hline
(\ref{6}.\ref{6-1}) \\
(\ref{6}.\ref{6-2}) \\
(\ref{6}.\ref{6-3}) \\
\hline
(\ref{5-2}.\ref{5-1-1}) \\
(\ref{5-2}.\ref{5-1-2}) \\
\hline
(\ref{5-4}.\ref{5-2-1}) \\
(\ref{5-4}.\ref{5-2-2}) \\
\hline
(\ref{4-1}) \\
\hline
(\ref{4-3}) \\
\hline
(\ref{4-2}) \\
\hline
(\ref{4-4}) \\
\hline
(\ref{3}) 
\end{array}
&
\begin{array}{c|c|c}
\deg(c_1) & c_1 & c_2 \\
\hline
\hline 
6 & 2H & 5 \\
6 & H+T & 4 \\
6 & H+C+L & 3 \\
\hline
5 & H+C & 3 \\
5 & H+L_1+L_2 & 2 \\
\hline
1 & H-C & 1 \\
1 & H-L_1-L_2 & 0 \\
\hline
4 & H+L & 2 \\
\hline
3 & H & 2 \\
\hline
2 & C & 1 \\
\hline
0 & 0 & 1 \\
\hline
3 & T & 1
\end{array}
\end{array}
\end{equation}
\end{scriptsize}
for some $T \in \TX$, $C \in \CX$, $L,\,L_1,\,L_2\in \LX$, with $C \cdot L = L_1 \cdot L_2 = 0$.
\end{thm}

\begin{corol}[Extensibility to the cubic threefold] Let $\EE$ be as above.
Then $\EE$ extends to the general cubic threefold $Y$ if and only if it is of type (\ref{4-4}), (\ref{4-3}), (\ref{6}.\ref{6-1}).
\end{corol}

\begin{proof}
By a result of Arrondo and Costa, \cite{arrondo-costa}, there are precisely $3$ families of (normalized)
indecomposable rank $2$ aCM bundles on the {\em general} cubic threefold $Y$,
corresponding to a line, a conic, and a linearly normal elliptic quintic contained in $Y$.
Their restriction to a general hyperplane section $X\subset Y$ corresponds respectively 
to bundles of type (\ref{4-4}), (\ref{4-3}), (\ref{6}.\ref{6-1}).

On the other hand, if $\EE$ belongs to a class other than the ones mentioned above, then $c_1(\EE)$
does not lift to a divisor class on $\p^3$. Since $\Pic(Y)$ is generated by the hyperplane class, $c_1(\EE)$
does not lift to $\Pic(Y)$, so that $\EE$ does not lift to $Y$ either.
\end{proof}

Supposing that a general section $s$ of $\EE$ vanishes in codimension $2$, we can write the following
exact sequences:
\begin{align}
\label{koszul-E} & 0 \rr \wedge^2(\EE^*) \rr \EE^* \rr J_Z \rr 0 \\
\label{ideal-Z}  & 0 \rr J_Z \rr \OO_X \rr \OO_Z \rr 0
\end{align}
where $Z\subset X$ is a subscheme of codimension $2$, with $c_2(\EE)=\len(Z) \geq 0$.
We first write a couple of lemmas.

\begin{lem} \label{its-extension}
Suppose that $\wedge^2(\EE) \simeq \LL_1 \ts \LL_2$, and let $\CC_2 \in |\LL_2|$ be a smooth rational curve
with $Z \subset \CC_2$.
Assume $\HH^1(\LL_1^*)=0$. Suppose further that $\len(Z) = \LL_1 \cdot \LL_2$.

Then we have an exact sequence:
\begin{equation} 
  0 \rr \LL_2^* \rr \EE^* \rr \LL_1^* \rr 0
\end{equation}
\end{lem}

\begin{proof}
Tensorize by $\LL_2$ the exact sequence (\ref{koszul-E}) and take global sections.
The hypothesis $\HH^1(\LL_1^*)=0$ implies that any nonzero morphism $\LL_2^* \rightarrow J_Z$ lifts to an
injective morphism $\LL_2^* \rightarrow \EE^*$.
On the other hand, since $Z$ is contained in $\CC_2$, we have $\LL_2^* \hookrightarrow J_Z$; denote by $\FF$ the cokernel of
$\LL_2^* \rr \EE^*$ and write the exact diagram (omitting zeroes all around):
\begin{equation} 
  \xymatrix{
                                        & \LL_2^* \ar@{=}[r] \ar[d] & \LL_2^* \ar[d] \\
  \LL_1^* \ts \LL_2^* \ar[r] \ar@{=}[d] & \EE^* \ar[r] \ar[d]       & J_Z \ar[d]     \\
  \LL_1^* \ts \LL_2^* \ar[r]            & \FF \ar[r]                & J_{Z,\CC_2}
  }
\end{equation}

Since $\CC_2$ is smooth and rational, by $\len(Z) = \LL_1 \cdot \LL_2$ we get
$J_{Z,\CC_2} \ts \LL_1 \simeq \OO_{\CC_2}$.
We have $\HH^1(\LL_2^*)=0$, so we can use Lemma \ref{trivial}, part (\ref{trivial-2}) to conclude that $\FF \ts \LL_1 \simeq \OO_X$
and thus recover the required exact sequence. 
\end{proof}

\begin{rmk} \label{its-extension-plus}
The result of the previous lemma holds also under the following weaker hypothesis:
\begin{itemize}
\item the curve $\CC_2$ is a connected union of smooth rational curves $\DD_j$;
\item the linear system $|\LL_1|$ contains a reduced connected curve $\CC_1$;
\item for each $j$ we assume $\DD_j \cdot \CC_1 = \len(Z_j)$ where $Z_j = Z \cap \DD_j$, and $Z = \cup_j Z_j$.
\end{itemize}
\end{rmk}

\begin{lem} \label{ctwoisone}
Let $\FF$ be a vector bundle on $X$, satisfying the following:
\begin{align}
& \HH^0(\wedge^2(\FF)(-1)) = 0 && \HH^1(\wedge^2(\FF)(-1)) = 0 \\
& 4 \leq \deg(c_1(\FF)) \leq 6 && c_2(\FF) = 1 \\
& \HH^1(\FF) = 0 && 
\end{align}
and assume that a general section of $\FF$ vanishes in codimension 2. Then $\FF$ is decomposable.
\end{lem}

\begin{proof}
Let $x\in X$ be the vanishing locus of a general global section $s$ of $\FF$.
Set $d=\deg(c_1(\FF))$ and $\MM = \wedge^2(\FF)$.
By Lemma \ref{classify-deg-5}, our hypothesis give $c_1(\MM) \in \RdX$.
So $\hh^0(\FF) = \hh^0(\MM) = d$ and $\hh^0(J_x \ts \MM) = d-1$.

Recall Lemma \ref{codimension}, and choose a general reduced reducible curve $\CC$ in $|\MM|$,
having an irreducible component $\CC_0$ containing $x$.
Let $\CC = \CC_0 \cup \DD$, and take $\CC'$ to be the union of $\CC_0$ and all but one
connected components of $\DD$ meeting $\CC_0$. This is possible for the components of $\CC$
form a simply connected graph. 
The divisor class $\CC'$ lies in ${\sf R}_{e}(X)$, for some $e < d$, and we denote it by $R_{e}$.

On the other hand, the class of the
remaining component of $\DD$ lies in ${\sf R}_{d-e}(X)$, and we denote it by $R_{d-e}$.
We have $R_{d} \cdot R_{d-e} =1$.

Clearly, in this process we can achieve $e \geq d-e$, i.e. $d \leq 2 \, e$.
Since all components of $\CC$ are smooth rational curves, we conclude by Remark
\ref{its-extension-plus} that $\FF^*$ fits into:
\[
0 \rr \OO_X(-R_e) \rr \FF^* \rr \OO_X(-R_{d-e}) \rr 0
\]

Now we look at the line bundle $\LL := \OO_X(R_e-R_{d-e})$:
here the bundle $\FF^*$ represents an element of $\HH^1(\LL^*)$.
We may assume $R_e \neq R_{d-e}$.
It is an easy exercise to prove:
\begin{align*}
& \chi(\LL^*) = d-e-2 && \chi(\LL) = e-2 \\
& \HH^0(\LL^*) = \HH^2(\LL^*) = 0  && \HH^2(\LL) = 0
\end{align*}

One sees easily that this implies our claim by Lemma \ref{other-splitting}, except in case
$d=4$, $e=2$. Notice that in this case we have $R_e-R_{d-e} = H-T$
(cfr. Remark \ref{combinatorics}: compose \ref{twisted-lines} with $\rho$) for some $T \in \TX$.
But $T$ is an aCM divisor, so we come to the same conclusion.
\end{proof}

From now on in this section, we let $\EE$ be an {\em indecomposable normalized rank $2$
aCM bundle} over $X$.

\begin{lem} \label{nowhere-vanishing}
Let $s$ be a nonzero global section of $\EE$.
\begin{enumerate}[i)]
\item \label{split} If $s$ is nowhere vanishing and $\hh^0(\EE)\geq 2$ then $\EE$ is decomposable.
\item \label{gg} If $\EE$ is globally generated, $\hh^0(\EE) \geq 4$ and $c_2(\EE)=1$, then $\EE$ is decomposable.
\end{enumerate}
\end{lem}

\begin{proof}
Statement (\ref{split}). We have an exact sequence:
\begin{equation} \label{nowhere}
0 \rr \wedge^2(\EE^*) \rr \EE^* \rr \OO_X\rr 0
\end{equation}

Since $\EE$ is aCM, this means that $\hh^1(\wedge^2(\EE^*))=1$.
On the other hand, since $\hh^0(\EE)\geq 2$, dualizing (\ref{nowhere})
we get that $\hh^0(\wedge^2(\EE)) \neq 0$.
This implies our statement by Lemma \ref{other-splitting}.

Statement (\ref{gg}) is an easy consequence of Lemma \ref{ctwoisone}
\end{proof}

\subsection{Linear resolutions}

We take into account the possible Chern classes of and aCM rank $2$ bundle $\EE$,
supposing that its minimal resolution is a $6 \times 6$ matrix of
linear forms.

\begin{prop} \label{chern-6}
Let $\EE$ have resolution (\ref{6}). Then one of the following cases must take place: 

\begin{scriptsize}
\begin{center}
$$\newcounter{cases-6}
\begin{array}{c|c|c|c|c}
    \mbox{Ref.} & c_1(\EE) & c_2(\EE) & \hh^0(\wedge^2(\EE(-H))) & \hh^1(\wedge^2(\EE)(-2\,H)) \\
    \hline 
    \hline
    \refstepcounter{cases-6} \label{6-1}
    \mbox{(A\arabic{cases-6})} & 2\,H & 5 & 4 & 0 \\
    \refstepcounter{cases-6} \label{6-2}
    \mbox{(A\arabic{cases-6})} & H+T  & 4 & 3 & 0 \\
    \refstepcounter{cases-6} \label{6-3}
    \mbox{(A\arabic{cases-6})} & H+C+L & 3 & 2 & 1
\end{array}$$
\end{center}
\end{scriptsize}
for some $T\in \TX$ $C\in \CX$, $L \in \LX$, with $C\cdot L=0$. We have:
 $$\fff(\EE) \,\, \text{skew-symmetric} \Longleftrightarrow c_2(\EE)=5
 \Longleftrightarrow c_1(\EE) \lequiv 2\, H$$ 
\end{prop}

\begin{proof}
Since $\EE$ has resolution (\ref{6}), it is globally generated.
Thus we can write down (\ref{ideal-Z}) and 
the Koszul sequence
(\ref{koszul-E}) associated to a section $s$ of $\EE$.
If $Z = \emptyset$, we conclude that $\EE$ is decomposable by Lemma \ref{nowhere-vanishing}.
So, we suppose that $Z$ consists of $\len(Z) = c_2(\EE)$, distinct points in $X$, with $c_2(\EE) \neq 0$.
By (\ref{ideal-Z}) and (\ref{koszul-E}) we have:
\begin{align}
\label{h1J} & \hh^0(J_Z)=\hh^2(J_Z)=0 && \hh^1(J_Z)=c_2(\EE)-1 \\
\label{h0L} & \hh^0(\wedge^2(\EE^*))=\hh^1(\wedge^2(\EE^*))=0 && \hh^2(\wedge^2(\EE^*))=\hh^1(J_Z)
\end{align}

Define the line bundle $\LL := \wedge^2(\EE)(-H)$.
Since $Z \neq \emptyset$ we have:
\begin{equation} \label{c2h0}
\hh^0(\LL) = c_2(\EE) - 1
\end{equation}

By Serre duality, $\LL$ satisfies the
hypothesis of Lemma \ref{classify-deg-5} and $\deg(\LL)=3$.
Therefore, assuming $\hh^0(\LL)\neq 0$, one of the alternatives (L\ref{LL-3-1}), \ldots, 
(L\ref{LL-3-4}) of this Lemma must take place.
The cases (L\ref{LL-3-1}), (L\ref{LL-3-2}), (L\ref{LL-3-3}) correspond to our table.
We must now exclude the possibilities:

\begin{enumerate}[I)]
\item \label{no-sections} No sections: $\hh^0(\LL)=0$, $c_2(\EE)=1$.
\item \label{one-section} One section: $\hh^0(\LL)=1$, $c_2(\EE)=2$ (i.e. case (L\ref{LL-3-4})).
\end{enumerate}

Case (\ref{no-sections}) is is obvious by Lemma \ref{ctwoisone}.
We prove Case (\ref{one-section}) by an {\em ad hoc} argument. Use the map \eqref{three-lines}
or Remark \ref{combinatorics} to see that $c_1(\EE) = T_1 + T_2$, for some $T_i \in \TX$
and $T_1 \cdot T_2 = 2$. For both $i$'s, since $\len(Z)=2$, the subscheme
$Z$ lies in a (reduced) curve $\CC_i$ of
$|\OO_X(T_i)|$. If $\CC_1$ and $\CC_2$ have no common component, the claim follows easily by
Remark \ref{its-extension-plus}.

So, assume that the $\CC_i$'s contain a common conic $\DD$ lying in $|\OO_X(C)|$, with $T_i = C+L^{(i)}$,
$L^{(i)} \in \LX$. Then, $L^{(1)} \cdot L^{(2)} = 0$.
Notice that a length $1$ subscheme of $Z$ is contained in $\DD$, while $Z$ itself can be contained or not in $\DD$.
In the former case, one can easily conclude by Remark \ref{its-extension-plus}.
In the latter case, we set
$\LL_1 = \OO_X(C)$, $\LL_2 = \OO_X(C+L^{(1)}+L^{(2)})$,
and our claim follows by the same remark.

On the other hand, assume that $\CC_1$ and
$\CC_2$, contain a line $L$ as a common component,
and set $C_i = T_i-L$. Then $L$ lies in $\{L_1,L_2,L_3\}$, with $T_1 + T_2 = H + L_1 + L_2 + L_3$.
So $Z$ cannot be contained  in $L$, for
we have $\HH^0(J_Z\ts \OO_X(T_1+T_2-H))=0$. By the Cayley-Bacharach property, no subscheme of $Z$ can
be contained in $L$, hence $Z \cap L = \emptyset$. But this is impossible since $C_1 \cdot C_2 = 1$.

The claim about skew-symmetry follows by Remark \ref{periodicity}.
\end{proof}

\subsection{Rank 5 resolution}

We will describe the behavior of the two sets of $5 \times 5$
resolutions, which are related by duality after Theorem \ref{full-classification}.

\begin{prop} \label{chern-5-1}
Let $\EE$ have resolution (\ref{5-2}). Then the possible cases are:
\begin{scriptsize}
\begin{center}
$$\newcounter{cases-5}
\begin{array}{c|c|c|c|c}
    \mbox{Ref.} & c_1(\EE) & c_2(\EE) & \hh^0(\wedge^2(\EE(-H))) & \hh^1(\wedge^2(\EE)(-2\,H)) \\
    \hline 
    \hline
    \refstepcounter{cases-5} \label{5-1-1}
    \mbox{(B.\arabic{cases-5})} & H+C & 3 & 2 & 0 \\
    \refstepcounter{cases-5} \label{5-1-2}
    \mbox{(B.\arabic{cases-5})} & H+L_1+L_2 & 2 & 1 & 1
\end{array}$$
\end{center}
\end{scriptsize}
for some $C,\,L_1,\,L_2 \subset X$, with $L_1 \cdot L_2 = 0$.
\end{prop}

\begin{proof}
The proof is analogous to that of Proposition \ref{chern-6}.
The exact sequences (\ref{koszul-E}), (\ref{ideal-Z}) and the formulas
(\ref{h1J}), (\ref{h0L}) still hold.
Set $\LL:= \wedge^2 \EE(-H)$. We still have the formula
(\ref{c2h0}). 
Again Lemma \ref{nowhere-vanishing} implies that $Z$ is nonempty,
so $c_2(\EE) \geq 1$.
On the other hand, the case $c_2(\EE)=1$ can be excluded by virtue of Lemma \ref{ctwoisone}.

So, the line bundle $\LL$ verifies the hypothesis of Lemma \ref{classify-deg-5},
with $\deg(\LL)=2$, which amounts to our table.
\end{proof}

\begin{prop} \label{chern-5-2}
Let $\EE$ have resolution (\ref{5-4}). Then the possible cases are:
\begin{scriptsize}
\begin{center}
$$\newcounter{cases-5-2}
\begin{array}{c|c|c|c|c}
    \mbox{Ref.} & c_1(\EE) & c_2(\EE) & \hh^0(\wedge^2(\EE)) & \hh^1(\wedge^2(\EE)(-H)) \\
    \hline 
    \hline
    \refstepcounter{cases-5-2} \label{5-2-1}
    \mbox{(C.\arabic{cases-5-2})} & L & 1 & 1 & 0 \\
    \refstepcounter{cases-5-2} \label{5-2-2}
    \mbox{(C.\arabic{cases-5-2})} & C-L & 0 & 0 & 1
\end{array}$$
\end{center}
\end{scriptsize}
for some $C,L \subset X$, with $C \cdot L = 1$.

Moreover, set $f:=\fff(\EE)_{|X}$, $g:=f^{\top}$, and $\FF :=
\cok(g)\ts \OO_X(-2)$.
Then $\FF$ is described by Proposition \ref{chern-5-1}, 
and we have:
\begin{equation} \label{E-dual-F}
\EE \simeq \FF^*(H) \qquad c_2(\FF)=c_2(\EE)+2
\end{equation}

Finally, there is minimal $2$-periodic exact sequence:
$$ \cdots \xr{f} \OO_X(-3) \oplus \OO_X(-4)^4 \xr{g} \OO_X(-2)^5
\xr{f} \OO_X \oplus \OO_X(-1)^4 \xr{g} \cdots $$
\end{prop}

\begin{proof}
Recall from Theorem \ref{full-classification} that the minimal graded
free resolution of the normalization of $\EE^*$ takes the form
(\ref{5-2}).
Grothendieck duality implies that the resolution matrix of $\EE^*$ is $\fff(\EE)^{\top}$.

Therefore we have (\ref{E-dual-F}).
This implies minimality in the $2$-periodic exact sequence provided by
Remark \ref{periodicity}, so we have the required exact infinite
which can be also written as:
$$0 \rr \EE^*(-2) \rr \OO_X(-2) \oplus \OO_X(-1)^4
 \xr{\fff(\FF)_{|X}=f^{\top}} \OO_X^5 \rr \FF \rr 0$$ 

The Chern classes of $\EE$ are thus computed in terms of those of
$\FF$ after Proposition \ref{chern-5-1}.
\end{proof}

\subsection{Rank 4 resolution}

In the next proposition we determine the Chern classes of
a bundle of type (\ref{4-3}) or (\ref{4-4}).
Denote it accordingly by $\EE$ or $\GG$.
Notice that if the matrix $\fff(\EE)$ is skew-symmetric, then
we obviously have $\gnrt(\EE)^* \simeq \syz(\EE)(t)$ for some $t$.
A consequence of the next proposition is that, in case $\rk(\gnrt(\EE))=4$,
this condition is also sufficient.

\begin{prop} \label{orientable-4}
Let $\EE$ and $\GG$ be as above.
\begin{enumerate}
\item \label{4-2-points} In case (\ref{4-3}) we have $c_1(\EE)\lequiv H$, $c_2(\EE)=2$.
\item \label{4-point} In case (\ref{4-4}), we have $c_1(\GG) = 0$, $c_2(\GG)=1$.
\end{enumerate}

In both cases the resolution matrix is skew-symmetric.
\end{prop}

\begin{proof}
Here we need a slightly more complicated proof,
due to the fact that in neither case is the bundle globally generated.
Let us look at (\ref{4-point}), and take a section $\sigma$ of $\GG$ (unique up to a scalar).
Write $\sigma^{\top}$ as a map $\GG \rr \wedge^2(\GG)$.
Then $\im(\sigma^{\top}) \subset \wedge^2(\GG)$ is isomorphic to $J_W \ts \NN$
where $W$ is a subscheme of codimension $2$ in $X$ (possibly empty), and $\NN$ is a line bundle.
Setting $\MM := \ker (\sigma^{\top})$, we get $\NN \simeq \MM^* \ts \wedge^2(\GG)$.
We have:
\begin{equation}
0 \rr \MM \rr \GG \xr{\sigma^{\top}} J_W \ts \NN \rr 0 \qquad \len(W)=1+\frac{\NN^2 + \MM^2}{2}
\end{equation}
where the second formula follows from Riemman-Roch and $\chi(\GG)=1$.
Since $\deg(\wedge^2 \GG) = 0$, assuming $\HH^0(J_W \ts \NN) \neq 0$ we obtain $\NN \simeq \OO_X$
and $W = \emptyset$. Then applying Lemma \ref{other-splitting} to $\MM(1)$ we deduce that $\GG$ is decomposable.

Thus we can assume $\HH^0(J_W \ts \NN) = 0$, $\HH^1(\MM)=0$, $\hh^0(\MM)=1$,
and we would like to prove that one of following two cases occurs:
\begin{small}
\begin{enumerate}[i)]
\item $\MM \simeq \OO_X$, $\NN \simeq \OO_X$ $\len(W) = 1$, or
\item $\MM \simeq \OO_X(L)$ $\NN \simeq \OO_X(-L)$, $\len(W) = 0$, for some $L \in \LX$. 
\end{enumerate}
\end{small}

Notice that $\GG(1)$ is indeed globally generated, which implies that the degree $3$ line bundle
$\wedge^2(\GG)(1)$ satisfies the
hypothesis of Lemma \ref{classify-deg-5}. A general section of $\GG(1)$
provides a map $\GG \rr \wedge^2(\GG)(1)$, hence 
an injective map $\MM \rr \wedge^2(\GG)(1)$. Notice that this map cannot be an isomorphism,
for otherwise $\GG$ would split as $\MM \oplus \wedge^2(\GG)(1)$. Thus $\deg(\MM) \leq 2$.
So, unless $\MM\simeq \OO_X$, Lemma \ref{classify-deg-5} applies to $\MM$.
Notice also that we have $\hh^0(J_W \ts \NN(1)) = 7 - \hh^0(\MM(1))$, and $\HH^1(\NN(1))=0$.
Summing up we are left with the cases:

\begin{equation} \label{cases-for-MM}
\begin{scriptsize}
\begin{tabular}{c|c|c|c|c}
$c_1(\MM)$ & $\deg(\MM)$ & $\deg(\NN(1))$ & $\hh^0(\MM(1))$ & $\hh^0(J_W \ts \NN(1))$ \\
\hline
\hline
$0$ & 0 & 3 & 4 & 3 \\
$L$ & 1 & 2 & 5 & 2 \\
$L_1 + L_2$ & 2 & 1 & 6 & 1 
\end{tabular}
\end{scriptsize}
\end{equation}

Suppose first $\MM \simeq \OO_X$. Here also $\NN(1)$ is classified by Lemma \ref{classify-deg-5},
so by $\hh^0(\NN(1))\geq 3$, we deduce that $c_1(\NN(1))$ is either $H$ or $T$, for some $T\in \TX$.
Accordingly,
$\NN^2$ equals $0$ or $-2$. In the former case we get $\len(W)=1$, which corresponds to our statement.
In the latter case, we get $\len(W)=0$, so $\GG$ is splits as $\OO_X \oplus \OO_X(T-H)$.

Looking at the remaining cases of (\ref{cases-for-MM}), we assume $c_1(\MM)=L$.
Then $\hh^0(\NN(1))\geq 2$ gives $c_1(\NN(1))=C$, for some conic $C$, so $c_1(\NN) = -L'$.
This implies $\len(W)=0$, and $\GG$ is decomposable unless $L' = L$.
We leave the last case to the reader.

\medskip

Finally, let us consider (\ref{4-2-points}). We define $\EE$ as the
normalization of the aCM bundle $\ker(\ppp(\GG)_{|X})$.
In case $c_1(\EE)\lequiv H$ we are done.
Excluding all other possibilities amounts to proving that $\hh^0(\wedge^2(\EE)) \leq 3$ leads
to a contradiction.
  
For instance, in case $\hh^0(\wedge^2(\GG))=3$, i.e. $c_1(\GG)\lequiv T$, we
know by the above discussion that $\GG$ splits as $\OO_X \oplus \OO_X(T-H)$ or
$\OO_X(L) \oplus \OO_X(C-H)$.
So the matrices $\fff(\GG)$ and $\fff(\EE)$ decompose in diagonal blocks.
This implies 
either that $\fff(\EE)$ can be reduced to the form (\ref{3}), either that
$\EE$ splits as $\OO_X(C) \oplus \OO_X(L)$. The remaining cases are similar.

The assertion about skew-symmetry follows from Remark \ref{periodicity}.
\end{proof}

In the next proposition, we examine the second pair of $4\times 4$
resolution matrices, i.e. (\ref{4-1}) and (\ref{4-2}) in Theorem
\ref{full-classification}.

\begin{prop} \label{chern-4-2}
Let $\EE$ and $\FF$ be an indecomposable rank $2$ aCM bundles on $X$
such that the minimal graded free resolution of $\EE$ (resp. of $\FF$) 
takes the form (\ref{4-1}) (resp. (\ref{4-2})).
Then we have:
\begin{align}
& c_1(\EE) \lequiv L+H && c_2(\EE)=2 \\
& c_1(\FF) \lequiv H-L && c_2(\FF)=1
\end{align}
for some $L \subset X$. Moreover $\EE$ and $\FF$ are related by:
\begin{align}
\label{E-dual-F-2} & \EE^* \simeq \FF(-1)  && \fff(\EE)^{\top}=\fff(\FF(2)) && \ker(\ppp(\EE)_{|X}) \simeq \FF(-1)
\end{align}
\end{prop}

\begin{proof}
Take a general section $s$ of the globally generated bundle $\EE$ and
write the Koszul resolution (\ref{koszul-E}).
Set $\LL:=\wedge^2(\EE)(-H)$.
By Lemma \ref{nowhere-vanishing}, the vanishing locus $Z$ in nonempty.
It follows by Lemma \ref{classify-deg-5} that $\hh^0(\LL) \neq 0$ if
and only if $c_1(\LL)\lequiv L$. By
the formula (\ref{c2h0}) we have $c_2(\EE)=\len(Z)=2$.

On the other hand, one excludes the case $c_2(\EE)=1$ by Lemma \ref{ctwoisone}.
All the remaining statements are clear by Remark \ref{periodicity},
and formulas (\ref{grothendieck}), (\ref{serre-duality}).
\end{proof}

\section{Extensions} \label{ext}

In this section we take into account rank $2$ aCM vector bundles on $X$
arising from extensions of aCM two line bundles $\MM$ and $\NN$.
These are also called {\em layered sheaves}.
In the following theorem, we show what is the resolution of such sheaf, according to $\MM$, $\NN$
and $\MM \cdot \NN$.

\begin{thm} \label{summary-extensions}
All cases of Theorem \ref{full-chern} are realized by an aCM indecomposable unobstructed rank $2$ bundle
$\EE$ on $X$ which is a nonsplitting extension:
\begin{equation} \label{the-extension}
0 \rr \MM \rr \EE \rr \NN \rr 0
\end{equation}
where $\MM$ and $\NN$ are aCM line bundles according to the following tables:

\begin{align}
\begin{tiny}
\label{normalized} \fbox{\txt{$\NN$ and $\MM$ \\ normalized}}
\end{tiny}
&&  
\begin{tiny}
\begin{tabular}{c||c|c|c|c|c|c|c|c}
Ref. & 	(\ref{6}.\ref{6-1},\ref{6-2},\ref{6-3}) & (\ref{5-2}.\ref{5-1-1}) & (\ref{4-1}) & (\ref{5-2}.\ref{5-1-1},\ref{5-1-2})
& (\ref{4-3}) & (\ref{4-2}) & (\ref{4-1}) & (\ref{4-3},\ref{3}) \\  
\hline
\hline
$c_1(\MM)$	& $T_1$ & $T$ & $C_1$ & $C$ & $C$ & $L_1$ & $L$ & $L$ \\
\hline
$c_1(\NN)$ & $T_2$ & $C$ & $C_2$ & $T$ & $L$& $L_2$ & $T$ & $C$ \\
\hline
$\MM \cdot \NN$ & (5,4,3) & 3 & 2 & (3,2) & 2 & 1 & 2 & (2,1)
\end{tabular}
\end{tiny} \\
\begin{tiny}
\label{not-normalized} \fbox{\txt{either $\NN$ or $\MM$ \\ not normalized}}
\end{tiny}
&&  
\begin{tiny}
\begin{tabular}{c|c|c|c|c|c|c}
Ref. & (\ref{4-3},\ref{3}) & ({\ref{4-2}}) & (\ref{5-4}.\ref{5-2-1},\ref{5-2-2}) & (\ref{4-4}) & (\ref{5-4}.\ref{5-2-1}) & (\ref{4-4}) \\
\hline
\hline
$c_1(\MM)$	& $T_1-H$ & $T-H$ & $T-H$ & $C-H$  & $L$ & $L$\\
\hline
$c_1(\NN)$ & $T_2$  & $C$ & $L$ & $L$ & $T-H$ & $C-H$ \\
\hline
$\MM \cdot \NN$ & (2,1) & 1 & (1,0) & 1 & 1 & 1
\end{tabular}
\end{tiny}
\end{align}

Moreover, these exhaust all indecomposable aCM rank 2 bundles which are extensions
of aCM line bundles.
\end{thm}

In spite of this theorem, we will see in the 
next section that {\em only some} indecomposable aCM rank 2 bundles are of this form.
Notice also that a given resolution can be obtained in different ways. 
The proof of the above result follows summarizing the propositions and remarks of this section.

Through this section we let $M = c_1(\MM)$, $N = c_1(\NN)$, $\Delta=M-N$ and $\Sigma= M+N$.

\begin{lem} \label{ext-gen}
Let $\MM$ and $\NN$ be line bundles on $X$, and suppose
$\Ext^1(\NN,\MM)= \HH^1(\MM \ts \NN^*)\neq 0$.
Let $\EE$ be a vector bundle corresponding to a nonzero element $[\EE] \in \Ext^1(\NN,\MM)$.
Then:
\begin{enumerate}[i)]
\item \label{ext-gen-1} 
we have $\HH^0(\MM\ts \NN^*) \simeq \HH^0(\EE^* \ts \MM) \simeq \HH^0(\EE \ts \NN^*)$;
\item \label{ext-gen-2} 
if $\HH^0(\MM\ts \NN^*)=\HH^0(\NN \ts \MM^*)=0$, then $\EE$ is
  simple;
\item \label{ext-gen-3} 
if $\HH^2(\MM \ts \NN^*)=\HH^2(\NN \ts \MM^*)=0$, then $\HH^2(\EEnd (\EE))=0$.
\end{enumerate}
\end{lem}

\begin{proof}
The bundle $\EE$ corresponds to an exact sequence:
\begin{equation} \label{M-N}
 0 \rr \MM \xr{\iota} \EE \xr{\sigma} \NN \rr 0
\end{equation}

So we have a commutative exact diagram (we omit zeroes all around the
diagram):
\begin{equation}
\xymatrix{
{\MM \ts \NN^*} \ar[r]^-{\iota} \ar[d]^-{\sigma^{\top}} & {\EE \ts \NN^*} \ar[r]^-{\sigma} \ar[d]^-{\sigma^{\top}} & {\OO_X} \ar[d]^-{\sigma^{\top}} \\
{\EE^* \ts \MM} \ar[r]^-{\iota} \ar[d]^-{\iota^{\top}} & {\EEnd{\EE}} \ar[r]^-{\sigma} \ar[d]^-{\iota^{\top}} & {\EE^*\ts
  \NN} \ar[d]^-{\iota^{\top}} \\
{\OO_X} \ar[r]^-{\iota} & {\EE \ts \MM^*} \ar[r]^-{\sigma} & {\NN \ts \MM^*}
}
\end{equation}

Since $0 \neq [\EE] \in \Ext^1(\NN,\MM)$, the boundary map
$\kk = \HH^0(\OO_X) \rr \HH^1(\MM \ts \NN^*)$ is nonzero.
So, taking cohomology of the left column and of the top row we obtain
(\ref{ext-gen-1}), which in turn implies (\ref{ext-gen-2}).
Since $\hh^1(\OO_X)=\hh^2(\OO_X)=0$, we also get (\ref{ext-gen-3}).
\end{proof}

\subsection{Pairs of twisted cubics}

Let $T_1,T_2 \in \TX$ and assume $M=T_1, N=T_2$.
If $T_1 \cdot T_2 \in \{1,2\}$, there are no extensions to examine
by Lemma \ref{all-extensions}. This proposition takes into account the remaining cases.

\begin{prop} \label{cubics-ext}
Let $\ell \in \{3,4,5\}$ and suppose $T_1 \cdot T_2 = \ell$. 
Then any nonsplitting extension $\EE$
between $\OO_X(T_1)$ and $\OO_X(T_2)$ is simple and unobstructed of type (\ref{6}.$6-\ell$).
The matrix $\fff(\EE)$
is skew-symmetric if and only if $T_1 \cdot T_2=5$.
\end{prop}

\begin{proof}
Clearly $T_1-T_2$ and $T_2-T_1$ are not effective unless $T_1 = T_2$, but in this case $T_1\cdot T_2 = 1$.
So Lemma \ref{ext-gen}
asserts that $\EE$ is simple and unobstructed.
We find
  $\hh^0(\EE(-H))=0$ and $\hh^0(\EE)=6$.
So, since is $\EE$ aCM, the resolution takes the required form.
Finally, $\fff(\EE)$ is skew-symmetric if and only if the line
bundle $\OO_X(T_1+T_2) \simeq \wedge^2 \EE$ lifts to a line bundle on $\p^3$.
This happens if and only if $T_1 + T_2 = 2\,H$, i.e. if and
only if $T_1 \cdot T_2=5$.
\end{proof}

Notice that, when $T_1 \cdot T_2 =5$,
we have $\fff(\OO_X(T_1))^{\top}=\fff(\OO_X(T_2 + H))$.
So, in this case, even if $\EE$ splits as $\OO_X(T_1) \oplus
\OO_X(T_2)$, the matrix $\fff(\EE)$ is skew-symmetric.

\begin{rmk}
Let $T_1$, $T_2$ be twisted cubics in $X$ with $T_1 \cdot T_2 \geq 3$,
and let $\EE$ correspond to a nonzero element
$[\EE] \in \Ext^1(\OO_X(T_2),\OO_X(T_1-H))$.
Then we have:
\begin{align*}
T_1 \cdot T_2 = 5 & \Longrightarrow c_1(\EE) = H  \Longrightarrow  \mbox{$\EE$ of type (\ref{4-3})}; \\
T_1 \cdot T_2 = 4 & \Longrightarrow c_1(\EE) = T_3 \Longrightarrow  \mbox{$\EE$ of type (\ref{3})}; \\
T_1 \cdot T_2 = 3 & \Longrightarrow c_1(\EE) = L+C  \Longrightarrow  \mbox{$\EE$ splits as $\OO_X(L) \oplus \OO_X(C)$}.
\end{align*}
\end{rmk}

\subsection{Pairs twisted cubic vs line, or twisted cubic vs conic}

We have seen in the previous subsection the first column of the
tables (\ref{normalized}) and (\ref{not-normalized}) contained in Theorem \ref{summary-extensions}.
Here we let $M-t\,H \in \TX$, for some $t\in \Z$.

\begin{prop} \label{extensions-T-T}
Let $\EE$ be a nonsplitting extension as (\ref{the-extension}), with $M=T+t\,H$
for some $t\in \Z, T \in \TX$, and
$N=\CX \cup \LX$.
Then $\EE$ is unobstructed and behaves according to the table:
\begin{center}
\begin{tiny}
\begin{tabular}{c|c|c|c|c|c|c|c}
$N$ & $t$ & $ T \cdot N$ & $\hh^1$ & $c_1(\EE)$ & $c_2(\EE)$ & \mbox{\rm type of $\EE$}
	& $\mbox{\rm simple}$ \\
\hline
\hline
$C_1$ & $0$ 	& $3$ & $1$ & $C_2+H$ & $3$ & (\ref{5-2}.\ref{5-1-1}) & $\checkmark$ \\
$C_1$ & $-1$	& $3$ & $2$ & $C_2$ & $1$ & (\ref{4-2}) & $\checkmark$ \\
$C$ 	& $-1$	& $2$ & $1$ & $L_1+L_2$ & $0$ & $\OO_X(L_1) \oplus \OO_X(L_2)$ &  \\
\hline
$L_1$ & $-1$ & $2$ & $2$ & $L_2$ & $1$ & (\ref{5-4}.\ref{5-2-1}) & $\checkmark$ \\
$L$ & $-1$ & $1$ & $1$ & $-L_1-L_2+H$ & $0$ & (\ref{5-4}.\ref{5-2-2}) & $\checkmark$ \\
$L$ & $-2$ & $2$ & $1$ & $-C$ & $0$ & $\OO_X \oplus \OO_X(-C)$ & 
\end{tabular}
\end{tiny}
\end{center}
for some $C_i\in \CX$, $L_j \in \LX$, with $L_1 \cdot L_2 = 0$, and 
where the column $\hh^1$ indicates the rank of $\HH^1(\MM\ts \NN^*(t))$, and {\em type}
shows the minimal resolution of $\EE$, or its splitting type.
\end{prop}

\begin{proof}
By the computations of Lemma \ref{all-extensions},
it is easy to check the data contained in the following table:
\begin{center}
\begin{tiny}
\begin{tabular}{c|c|c|c|c|c|c|c}
$N$ & $t$ & $T \cdot N$ & $\hh^1$ & $T-N+H$ & $(T-N+H)^2$ & $c_1(\EE)$ & $c_2(\EE)$ \\
\hline
\hline
$C_1$ & $0$ & $3$ & $1$ & $T+L$ & $-1$ & $C_2+H$ 	& $3$ \\
$C_1$ & $-1$ & $3$ & $2$ & $T+L$ & $-3$ & $C_2$ 		& $1$ \\
$C$ 	& $-1$ & $2$ & $1$ & $R_4$ & $-1$ & $L_1+L_2$& $0$ \\
\hline
$L_1$ & $-1$ & $2$ & $2$ & $L_2$ 			& $-1$ & $L_2$ 				& $1$ \\
$L_1$ & $-1$ & $1$ & $1$ & $H-L_3-L_4$& $-3$ & $H-L_2-L_3$ 	& $0$ \\
$L_1$ & $-2$ & $2$ & $1$ & $L_2$ 			& $-1$ & $L_2$ 				& $0$ \\
\end{tabular}
\end{tiny}
\end{center}

The column $T-N+H$ gives the value of $\Delta$, and $c_1(\EE)$ is derived
from $\Sigma$ of Lemma \ref{all-extensions}.
In view of this, simplicity and unobstructedness
follow immediately in the required cases by Lemma \ref{ext-gen}.
Further, $\EE$ has the desired minimal resolution by a Hilbert polynomial computation.
For the splitting type of $\EE$, apply Lemma \ref{other-splitting} (\ref{two-curves}) to $\MM\ts \NN^*$
to recover the splitting type of $\EE \ts \NN^*$.
Use then the maps of Remark \ref{combinatorics} to
write down the summands of $\EE$ (namely we need \eqref{twisted-lines}, $\rho$ and $\tau$).
\end{proof}

Analogously, we may assume $N \in \TX$.
We leave the proof of the following proposition as an exercise,
noting that the splitting type is given by Lemma \ref{other-splitting},
and the maps \eqref{C+L=T}, \eqref{L+C=T} of Remark \ref{combinatorics}.

\begin{prop}
Let $\EE$ be a nonsplitting extension as (\ref{the-extension}), with $M-t\,H \in \CX \cup \LX$
for some $t\in \Z$ and $N= T \in \TX$. Set $M'=M-t\,H$.
Then $\EE$ is unobstructed and behaves according to the table:
\begin{center}
\begin{tiny}
\begin{tabular}{c|c|c|c|c|c|c|c}
$M'$ & $t$ & $ M' \cdot N$ & $\hh^1$ & $c_1(\EE)$ & $c_2(\EE)$ & \mbox{\rm type of $\EE$}
	& $\mbox{\rm simple}$ \\
\hline
\hline
$C_1$ & $0$ 		& $3$ & $2$ & $C_2+H$ 		& $3$ & (\ref{5-2}.\ref{5-1-1}) & $\checkmark$ \\
$C_1$ & $-1$ 		& $3$ & $1$ & $C_2$ 			& $0$ & $\OO_X \oplus \OO_X(C_2)$ &  \\
$C$ 	& $0$ 		& $2$ & $1$ & $L_1+L_2+H$ & $2$ & (\ref{5-2}.\ref{5-1-2}) &  $\checkmark$\\
\hline
$L_1$ &	$0$ & $2$ & $2$ & $L_2+H$ & $2$ & (\ref{4-1}) & $\checkmark$ \\
$L_1$ &	$1$ & $2$ & $1$ & $L_2$ 	& $1$ & (\ref{5-4}.\ref{5-2-1}) & $\checkmark$ \\
$L$ &		$0$ & $1$ & $1$ & $C_1+C_2$ & $1$ & $\OO_X(C_1) \oplus \OO_X(C_2)$ & 
\end{tabular}
\end{tiny}
\end{center}
for some $C_i\in \CX$, $L_i \in \LX$, with $L_1 \cdot L_2 = 0$, $C_1\cdot C_2=1$,  
where the column $\hh^1$ indicates the rank of $\HH^1(\MM\ts \NN^*(t))$, and {\em type}
shows the minimal resolution of $\EE$, or its splitting type.
\footnote{Here we implicitly normalize $\EE$ in the case with $t=1$, giving rise to (\ref{5-4}.\ref{5-2-1}).}
\end{prop}

\subsection{Lines and conics}

In the following proposition we analyze extensions of pairs of line bundles in $\LX$ or $\CX$.
The proof is analogous to that of \ref{extensions-T-T}.

\begin{prop}
Let $\EE$ be a nonsplitting extension as (\ref{the-extension}), with $M-t\,H$ and $C \in \LX$ or
$M-t\,H$ and $N \in \CX$, for some $t\in \Z$. Set $M'=M-t\,H$.
Then $\EE$ is unobstructed and described by the table:
\begin{center}
\begin{tiny}
\begin{tabular}{c|c|c|c|c|c|c|c}
$N$ & $t$ & $ M' \cdot N$ & $\hh^1$ & $c_1(\EE)$ & $c_2(\EE)$ & \mbox{\rm type of $\EE$}
	& $\mbox{\rm simple}$ \\
\hline
\hline
$C$ & $0$ 	& $2$ & $1$ & $L+H$ & $2$ & (\ref{4-1}) & $\checkmark$ \\
$C$ & $-1$	& $2$ & $1$ & $L$ 	& $0$ & $\OO_X \oplus \OO_X(L)$ & \\
\hline
$L_1$ & $0$ 	& $1$ & $1$ & $C$ 	& $1$ & (\ref{4-2}) & $\checkmark$ \\
$L_1$ & $-1$ 	& $1$ & $1$ & $-L_2$ 	& $0$ & $\OO_X \oplus \OO_X(-L_2)$ &
\end{tabular}
\end{tiny}
\end{center}
for some $L,L_2 \in \LX$, $C\in \CX$.
\end{prop}

\begin{rmk}
Taking $L_1,L_2$ in $\LX$ and $C_i = \rho(L_i)$, we have a natural isomorphism
 $$\psi : \HH^1(\OO_X(C_1-C_2)) \xr{\,\, \simeq \,\,}  \HH^1(\OO_X(L_2-L_1))$$

If the bundle $\EE$ corresponds to an element $[\EE]$ of
$\HH^1(\OO_X(L_2-L_1))$, then the bundle $\FF$ corresponding to $\psi([\EE])$
is isomorphic to $\EE^*(-H)$. Under this correspondence we have
$\fff(\EE)=\fff(\FF)^{\top}$ and there is a $2$-periodic minimal exact
sequence:

$$ \cdots \rr
\OO_X(-3)^4
\xr{\fff(\EE)_{|X}}
\begin{array}{c}
\OO_X(-2)^2 \\
\oplus \\
\OO_X(-1)^2
\end{array}
 \xr{\fff(\FF)_{|X}}
\OO_X^4
 \rr \FF \rr 0
$$
\end{rmk}

\begin{prop}
Let $\EE$ be a nonsplitting extension of the form (\ref{the-extension}). Assume
$M-t\,H \in \LX$ and $\rho(N)\in \LX$, or $M-t\,H \in \CX$ and $\rho(N)\in \CX$,
for some $t\in \Z$. Set $M'=M-t\,H$.
Then $\EE$ is unobstructed and we have the table:
\begin{center}
\begin{tiny}
\begin{tabular}{c|c|c|c|c|c|c|c|c}
$M'$ &$N$ & $t$ & $ M' \cdot N$ & $\hh^1$ & $c_1(\EE)$ & $c_2(\EE)$ & \mbox{\rm type of $\EE$}
	& $\mbox{\rm simple}$ \\
\hline
\hline
$C$ & $L$ & $0$ 	& $2$ & $1$ & $H$ 	& $2$ & (\ref{4-3}) & $\checkmark$ \\
$C$ & $L$ & $-1$	& $2$ & $2$ & $0$ 	& $1$ & (\ref{4-4}) & \\
$C$ & $L$ & $-1$ 	& $1$ & $1$ & $T-H$ & $2$ & $\OO_X \oplus \OO_X(T-H)$ &  \\
$C$ & $L$ & $-2$ 	& $2$ & $1$ & $-H$ 	& $0$ & $\OO_X \oplus \OO_X(-H)$ &  \\
\hline
$L$ & $C$ & $1$ 	& $2$ & $1$ & $2\,H$ 	& $4$ & (\ref{4-4}) &  \\
$L$ & $C$ & $0$ 	& $2$ & $2$ & $H$ 		& $2$ & (\ref{4-3}) & $\checkmark$ \\
$L$ & $C$ & $0$ 	& $1$ & $1$ & $T$ 		& $1$ & (\ref{3}) & $\checkmark$ \\
$L$ & $C$ & $-1$ 	& $2$ & $1$ & $0$ 		& $0$ & $\OO_X \oplus \OO_X$ & 
\end{tabular}
\end{tiny}
\end{center}
for some $L,L_2 \in \LX$, $C\in \CX$, $T \in \TX$.
\end{prop}

\begin{proof}
Using the classification of $\Delta$ and $\Sigma$ of Lemma \ref{all-extensions} one
can show easily all the information contained in our table, except for the indecomposability of
$\EE$ in case (\ref{4-4}).

So assume $M=-L$, $N=L\in \LX$, and suppose that $\EE$ decomposes as
a direct sum of two line bundles $\LL_1$ and $\LL_2$.
Notice that the $\LL_i$'s are aCM, and $\LL_2 \simeq \LL_1^*$.
So by Proposition \ref{acm-line-bundles} we can write: 
 $$c_2(\EE) = c_2(\LL_1 \oplus \LL_1^*) = - \LL_1^* = 2-m+3\,n^2$$
for some $n \in \Z$ and $m \in \{0,1,2,3\}$. Setting this number equal to $1$
we obtain $m=1$ and $n=0$ i.e. $\LL_1 \simeq \OO_X(L')$, for some
line $L' \in \LX$.
Then $\EE$ is indecomposable since $\hh^0(\EEnd (\OO_X(L')\oplus \OO_X(-L'))) = 3$,
while $\hh^0(\EEnd(\EE))=2$.

\end{proof}

\section{Moduli spaces} \label{moduli}

Here we draw a few remarks on moduli spaces of aCM
bundles. We only aim at some birational description of these families.
Through this section we let $\EE$ be an indecomposable rank $2$ aCM bundle on $X$.

\subsection{Moduli of linear resolutions}

We consider here moduli spaces of bundles whose
minimal graded free resolution is a $6 \times 6$ square matrix of
linear forms (type \ref{6}).

\begin{prop} \label{semistable-A}
Let $\EE$ be indecomposable of type (\ref{6}). Then $\EE$ is semistable.
\end{prop}

\begin{proof}
Since $\EE$ is globally generated we can write down (\ref{koszul-E}), where $Z$ consists of a set
of distinct points of $X$ of cardinality $c_2(\EE)$, satisfying the Cayley-Bacharach property for
the line bundle $\LL := \wedge^2(\EE)(-H)$.
Our assumption implies $\HH^0(J_Z\ts \LL) = 0$.
So, given any subscheme $W$ of $Z$ with $\len(W)=c_2(\EE)-1$,
we must have $\HH^0(J_W\ts \LL)=0$.

Consider now a destabilizing rank $1$ subbundle $\KK^*$ of $\EE^*$, with $\deg(\KK^*)\geq -2$.
Then $c_1(\KK)$ is an effective divisor class of degree at most $2$, and a curve in $|\KK|$ should contain $Z$.  
Now we look at the cases of Proposition \ref{chern-6} separately.

{\bf Case (\ref{6}.\ref{6-1})}.
We have $\HH^0(J_Z(H))= 0$, so $c_1(\KK) = L_1 + L_2$, $L_i \in \LX$, $L_1 \cdot L_2 = 0$.
At least $3$ points of $Z$ must lie on one of the lines $L_1$, $L_2$. So we find a subscheme $W$ of $Z$ made
of $4$ points in a plane. But this contradicts the Cayley-Bacharach property.

{\bf Case (\ref{6}.\ref{6-2})}. Again we have $\HH^0(J_Z(H))= 0$ and $c_1(\KK) = L_1 + L_2$ as above.
Let ${\sf U} \subset X$ be the union of the lines $L_1$ and $L_2$. Then we have the following exact diagram
(where we omit zeroes all around):
\begin{equation} \label{no-two-lines}
\xymatrix{
	& \OO_X(-{\sf U}) \ar@{=}[r] \ar[d] & \OO_X(-{\sf U}) \ar[d] \\
	\wedge^2(\EE^*) \ar@{=}[d] \ar[r] & \EE \ar[d] \ar[r] & J_Z \ar[d]\\ 
	\wedge^2(\EE^*) \ar[r] & \GG \ar[r] & J_{Z,{\sf U}}
	}
\end{equation}
for some coherent sheaf $\GG$. It suffices to show that the group $\HH^1(J_{Z,{\sf U}} \ts \LL)$ vanishes
under our assumptions. Indeed in this case the bottom row of (\ref{no-two-lines})
splits, which is not possible since we have $\Hom(\EE,\wedge^2(\EE^*))=0$. We may write
$J_{Z,{\sf U}}$ as $\OO_{L_1}(-Z_1) \oplus \OO_{L_2}(-Z_2)$, where $Z_i = Z \cap L_i$.

Notice that we must have $\len(Z_i)=2$ for each $i$,
for otherwise $Z$ is contained in a hyperplane.
Thus it suffices to prove that $\LL$ meets each $L_i$ at least at one point.

Recall that in our case $c_1(\LL)=T \in \TX$. If $T \cdot L_1 = 0$, then there exists a reducible curve
in $|\OO_X(T)|$ made of the union of $L_1$ and a conic though one point of $Z_2$. Thus the Cayley-Bacharach property is
violated and we are done.

{\bf Case (\ref{6}.\ref{6-3})}. Here we have $\hh^0(J_Z(H))=1$.
Consider a destabilizing rank $1$ subbundle $\KK^*$
of $\EE^*$. This time we must exclude the two cases $c_1(\KK)=L_1+L_2$ or $c_1(\KK) \in \CX$.
The former is similar to the one discussed above, so we omit it.
We are left with the latter: we set $c_1(\KK)=C_1$, and write $\CC_1$
for the curve in $|\OO_X(C_1)|$ containing $Z$.
Recall that here we have $c_1(\LL)=L + C$, with $L \in \TX$, $C \in \CX$.
Of course $C_1 \cdot C \neq 0$, for otherwise $Z$ is contained in a curve of $|\OO_X(C)|$.
Notice also that if $C_1 \cdot (C + L) \geq 2$, after writing down a diagram similar to
(\ref{no-two-lines}), with ${\sf U}$ replaced by a curve in $|\OO_X(C_1)|$,
we conclude that $\EE$ is decomposable.
So we must only
exclude the case $C_1 \cdot C = 1$, $C_1 \cdot L = 0$. Take the divisor $R_4 := H+L+C-C_1$,
observe that it lies in ${\sf R}_4(X)$ and contains $Z$.
Then, if $\CC_1$ is smooth, we can apply Lemma \ref{its-extension} and see that $\EE$ is decomposable.
To conclude the proof, 
assume that $\CC_1$ consists of two lines $L_3$ and $L_4$, and write $Z_i = Z \cap L_i$, $i\in \{3,4\}$.
Again by the Cayley-Bacharach property we easily see:

\begin{small}
\begin{align*}
& L \cdot L_3 =0 && L \cdot L_4 = 0 \\
& R_4 \cdot L_3 = \len(Z_3) && R_4 \cdot L_4 = \len(Z_4)
\end{align*}
\end{small}

Hence we are done by Remark \ref{its-extension-plus} and Lemma \ref{other-splitting}.
\end{proof}

By the previous Lemma, we can view a bundle $\EE$ of type (\ref{6}) as an element of 
$\MCMss(2;c_1(\EE),c_2(\EE))$, and we can consider the general element of this moduli space.

\begin{thm} \label{moduli-T-T}
In the three cases (\ref{6}.\ref{6-1}), (\ref{6}.\ref{6-2}), (\ref{6}.\ref{6-3}),
the general bundle $\EE$ is stable, $\MCMst(2;c_1(E),c_2(\EE))$
is a smooth irreducible rational variety of dimension $2\, c_2(\EE) - 5$.
\end{thm}

The proof is subdivided into the following Lemmas.

\begin{lem} \label{CB-T-T}
Let $\ell \in \{3,4,5\}$, $T \in \TX$, $L\in \LX$, $C\in \CX$, with $C \cdot L= 0$.
Set $\LL_5 = \OO_X(H)$, $\LL_4 = \OO_X(T)$, $\LL_3 = \OO_X(C+L)$.
Then there are open subsets $\mathcal{H}_{\ell} \subset \Hilb_{\ell}(X)$ such that, for any subscheme
$Z \in \mathcal{H}_{\ell}$, there exist a rank 2 aCM vector bundle $\EE$ on $X$
with $c_1(\EE)=c_1(\LL_{\ell})+H$, $c_2(\EE)=\ell$, and a section $s\in \HH^0(\EE)$
with $Z=\{ s=0 \}$.
\end{lem}

\begin{proof}
We define the open subsets $\mathcal{H}_{\ell}$ as follows.
\begin{small}
\begin{align*}
\widetilde{\mathcal{H}}_{\ell} & = \{Z \in \Hilb_{\ell}(X) \, | \,
	\mbox{$\forall \, W \subset Z$ with $\len(W)=\ell-1$, 
	we have $\HH^0(J_W\ts \LL_d)=0$} \} \\
\mathcal{H}_{5} & = \{Z \in \widetilde{\mathcal{H}}_{5} \, | \, 
	\mbox{$Z$ is contained in no hyperplane} \} \\
\mathcal{H}_{4} & = \{Z \in \widetilde{\mathcal{H}}_{4} \, | \, 
	\mbox{$Z$ is contained in no hyperplane} \} \\
\mathcal{H}_{3} & = \{Z \in \widetilde{\mathcal{H}}_{3} \, | \, 
	\mbox{$Z$ is contained in only one hyperplane} \}
\end{align*}
\end{small}

For $Z \in \mathcal{H}_{\ell}$ one easily checks  the Cayley-Bacharach property for the pair $(\LL,Z)$,
so $\EE$ is given by Theorem \ref{hartshorne-serre} and we have
the exact sequence (\ref{CB}).
It is easy to see that, for $Z \in \mathcal{H}_{\ell}$, we have:
  \begin{align}
  \label{single-extension}
  & \hh^1(J_Z \ts \LL_{\ell}) = 1 \\
  & \HH^1(J_Z (t)) = 0 && \mbox{for $t \geq 2$} \\
  & \HH^1(J_Z \ts \LL_{\ell}(t)) = 0 && \mbox{for $t \geq 1$}
  \end{align}

Therefore, since $\hh^1(\EE(t))=\hh^1(\EE^*(-t-1))$, in order to prove that $\EE$ is
aCM it suffices to show: 
  \begin{align} 
    \label{CB-acm-1} & \hh^1(\EE^*(1)) = 0 \\
    \label{CB-acm-2} & \hh^1(\EE(-1)) =0
  \end{align}

Condition (\ref{CB-acm-2}) holds whenever (\ref{single-extension}) holds.
Indeed, taking global section in (\ref{CB-dual}), the statement follows by the commutative diagram:
\begin{equation} \label{zero-h1}
  \xymatrix{
  {\kk \simeq \HH^1(J_Z \ts \LL_{\ell})} \ar[r] \ar_-{\mbox{Serre}}[d] & {\HH^2(\OO_X(-1)) \simeq \kk} \ar^-{\mbox{Serre}}[d] \\
  {\kk \simeq \Ext^1(J_Z,\LL_{\ell})^*} \ar^-{\simeq}[r] & {\HH^0(\OO_X)^* \simeq \kk}  
  }
\end{equation}

The bottom map is an isomorphism for it corresponds to the extension
given by $\EE$, which is nontrivial. So let us prove (\ref{CB-acm-1}).
If $\ell=5$, since $c_1(\EE) = 2\,H$, (\ref{CB-acm-2})
implies (\ref{CB-acm-1}). If $\ell=4$,
we have $\hh^1(J_Z(1))=0$ so (\ref{CB-acm-1}) holds.
We are left with the case $\ell = 3$.

Here we have $\hh^0(J_Z(1))=1$.
Assume that the boundary map
$\partial: \HH^0(J_Z(1)) \rr \HH^1(\LL_3^*)$
is zero. 
Notice that the generator of $\HH^1(\LL_3^*)$ corresponds to the vector bundle
$\FF:=\OO_X(-L-H) \oplus \OO_X(-C-H)$.
This implies that $\Hom(\FF,\EE^*) \simeq \Hom(\FF,J_Z)$, since $\hh^1(\FF^* \ts \wedge^2(\EE))=0$.
So, whenever the morphism $\OO_X(-1) \rr J_Z$ is nonzero, we lift it to a nonzero morphism $\FF \rr \EE$.
So $\partial \neq 0$ and we are done.
\end{proof}

\begin{lem}
Fix hypothesis as in Lemma \ref{CB-T-T}.
Then the aCM bundle $\EE$ is of type (\ref{6}).
For general $Z \in \mathcal{H}_{\ell}$, the bundle $\EE$ is stable.
The moduli space $\MCMst(2;c_1(\LL_{\ell}),\ell)$ is smooth and
irreducible of dimension $2\, \ell-5$.
\end{lem}

\begin{proof}
It is straightforward to see that, for $Z \in \mathcal{H}_{\ell}$, we have:
  \[\hh^0(\EE(-1))=0 \qquad  \hh^0(\EE)=6 \]

In view of Theorem \ref{full-classification}, this implies at once that 
the resolution of $\EE$ takes the form (\ref{6}).
By the previous Lemma and Proposition \ref{semistable-A}, we have the
semistable aCM vector bundle $\EE$.
Consider now the open subsets of $\mathcal{H}_{\ell}$ defined by:
\begin{small}
\[
\mathcal{H}_{\ell}^{\circ} = \left\{ Z \in \mathcal{H}_{\ell} \, \biggl\lvert \, 
	\begin{array}{l}
	\text{$Z$ is contained in no divisor $\DD$ of degree $3$,} \\
	\text{except $\DD \in |\OO_X(H)|$ if $\ell  = 3$}
	\end{array} \right\}
\]
\end{small}

For $Z \in \mathcal{H}_{\ell}^{\circ}$, the vector bundle given by Lemma \ref{CB-T-T}
is stable.
Our discussion implies that the map $\xi$ or Remark \ref{hilbert-serre} is dominant
and in fact the restriction of $\zeta$ to $\mathcal{H}_{\ell}^{\circ}$ is a birational morphism.
Since $\Hilb_{\ell}(X)$ is irreducible, the same holds for $\FMCMst(2;c_1(\EE),c_2(\EE))$ and thus
for $\MCMst(2;c_1(\EE),c_2(\EE))$ since $\eta$ is
a surjective map. Rationality of our moduli space also follows.
Now, by $\dim(\eta^{-1}([\EE]))=5$, we conclude that
$\MCMst(2;c_1(\EE),c_2(\EE))$ has dimension: 
 $$\dim(\Hilb_{\ell}(X))-5=2\, \ell-5$$
\end{proof}

In case $c_2(\EE)=5$, the moduli space $\MCMst(2;2\,H,5)$ can also be
described as the quotient of the space of skew-symmetric $6 \times 6$
matrices with linear entries and with Pfaffian equal to $F$ by the
action of $\SL(6)$ acting by conjugation.

\subsection{Families of rank 5 matrices}

We will analyze here the family of bundles of type (\ref{5-2}),
and we separate the cases (\ref{5-2}.\ref{5-1-1}) and (\ref{5-2}.\ref{5-1-2}).

\begin{prop}
Let $\EE$ be of type (\ref{5-2}.\ref{5-1-1}). The $\EE$ is stable unless it is an extension
of the form:
\begin{equation} \label{T-C-ext}
0 \rr \OO_X(T) \rr \EE \rr \OO_X(C) \rr 0
\end{equation}
for some $T \in \TX$, $C \in \CX$, $T \cdot C = 3$.
\end{prop}

\begin{proof}
We follow the proof of Proposition \ref{semistable-A}, and we introduce the same notation.
Notice that in this case $\EE$ is stable if and only if it is semistable.
Again the destabilizing subbundle $\KK^*$ of $\EE^*$ must be effective of degree at most $2$,
and one proves that it cannot be a sum of two skew lines. So let us assume $c_1(\KK)=C_1 \in \CX$.
Here as well, it is enough to study the case $C_1 \cdot C = 1$.

So we set $T = H + C - C_1$. It is easy to check that $T$ is a divisor class in $\TX$, with $T \cdot C = 3$.
Using the argument introduced at the end of the proof of Proposition \ref{semistable-A},
we obtain an extension of the form (\ref{C-T-ext}) in this case.
The difference here is that $\EE$ need not decompose, for $Z$ need not lie in a curve in $|\OO_X(T)|$.
\end{proof}

\begin{thm} \label{moduli-T-C}
Take $C \in \CX$, and set $L = \rho(C)$. 
Then the moduli space $\MCMst(2;C+H,3)$ of bundles of type (\ref{5-2}.\ref{5-1-1})
is isomorphic to $X \setminus L$.
\end{thm}

\begin{proof}
Recall the duality between bundles of type (\ref{5-2}.\ref{5-1-1}) and (\ref{5-4}.\ref{5-2-1}),
defined by $\EE \mapsto \FF = \EE^*(-H)$. It provides an isomorphism:
\[
\MCMst(2;C+H,3) \longleftrightarrow \MCMst(2;\rho(C),1)
\]

So take a bundle $\FF$ of type $(\ref{5-4}.\ref{5-2-1})$, with $c_1(\FF) = L = \rho(C)$ and $c_2(\FF)=1$,
and consider a nonzero global section $s$ of $\FF$ as 
a map $s : \FF^* \rr \OO_X$. Recall that the section $s$ is unique up to nonzero scalar.
Since $\deg(c_1(\FF))=1$, the vanishing locus $Z= \{s=0\}$ must have codimension two in $X$
(i.e. it must be a single point $z$ of $X$).
Indeed if the image of $s$ contained an effective divisor lying in $|\KK|$,
then the line bundle $\KK$ would destabilize $\FF$.
This gives a map:

\[
\xi : \MCMst(2;L,1) \rr X \qquad \mbox{defined by $\FF \mapsto \, z = \{s=0\}$}
\]

On the other hand, for each point $z \in X$, by Theorem \ref{hartshorne-serre} we obtain a locally free
sheaf $\GG$, given as the unique extension:
\[
0 \rr \OO_X(-L) \rr \GG^* \rr J_{z} \rr 0
\]

It is easy to see that $\GG$ is a stable sheaf which is aCM whenever $z$ does not lie in $L$.
However if $z$ does lie in $L$ then $\hh^1(\GG)=1$. This defines again the map
$\zeta: X \setminus L \rr \MCMst(2;L,1)$, which is clearly an inverse to $\xi$. This proves our result.
\end{proof}

For bundles of type (\ref{5-2}.\ref{5-1-2}), a different phenomenon occurs.

\begin{prop}
Let $\EE$ be of type (\ref{5-2}.\ref{5-1-2}). Then there is a pair $(T,C) \in \TX \times \CX$
such that $\EE$ fits into the following extension:
\begin{equation} \label{C-T-ext}
0 \rr \OO_X(C) \rr \EE \rr \OO_X(T) \rr 0
\end{equation}

There are $5$ different pairs $(T,C)$ that express $\EE$ as an extension of this form.
\end{prop}

\begin{proof}
We have $L_1+L_2 = c_1(\EE)-H$, with $L_i \in \LX$, $L_1 \cdot L_2 = 0$.
Let us borrow the notation again from the proof of Proposition \ref{semistable-A}.
Taking a general section of $\EE$, its vanishing locus $Z$ will consist of two distinct points of $X$.
Notice that, given any element $T \in \TX$,
we can find a curve in $\CC$ in $|\OO_X(T)|$ containing $Z$, i.e. an injection
$j : \OO_X(-T) \hookrightarrow J_Z$.
In order to lift $j$ to $\EE^*$ we have to check $\HH^1(J_Z(L_1+L_2-T))=0$.
One computes $\chi(L_1+L_2-T)=-T\cdot(L_1+L_2)$, so we must have $T\cdot L_1 = T \cdot L_2 = 0$.
Then $T-L_1-L_2$ is an element $L$ of $\LX$.
So we have to choose $L$ among the $5$ lines meeting both $L_1$ and $L_2$.

Given such choice for $L$ (equivalently, for $T$), we set $C = L_1+L_2+H-T$.
It is easy to see that $C$ lies in $\CX$, and that $T \cdot C = 2$.
If the curve $\CC$ containing $Z$ is smooth, we get our statement by Lemma \ref{its-extension}.
On the other hand, assume that $\CC$ contains a line $L_3$ with $T \cdot L_3 = 0$,
so that the class $D$ of the residual curve $\DD$ lies in $T-L_3 \in \CX$.
We have $L_3 \cdot (L_1 + L_2)=L_3 \cdot C - 1$, so $L_3 \cdot C = 0$ implies $L_3 \in \{L_1,L_2\}$.
So the Cayley-Bacharach property means that if $Z \cap L_3$ is nonempty, then $L_3 \cdot C \geq 1$.
Notice that $L_3 \cdot C = 2$ implies $C=D$, a contradiction.

Assume that $\DD$ is smooth. If $Z \subset L_3$, we set $R_4 = C+H-L_3$.
Observe that $R_4$ sits in ${\sf R}_4(X)$, and $L_3 \cdot R_4 = 2$, so using
Lemma \ref{its-extension}, one can easily see that $\EE$ is decomposable.
If $\len(Z \cap L_3)=1$, then we obtain $C \cdot L_3=1$, and we have the extension (\ref{C-T-ext})
by Remark \ref{its-extension-plus}.
Finally if $Z \cap L_3 = \emptyset$, we set $T' = H + C - D$ whence $T'\cdot D = 2$.
Again we can use Remark \ref{its-extension-plus} and conclude that $\EE$ is decomposable,
by $\HH^1(\OO_X(T'-D))=0$.

We leave it to the reader to work out the case when $\DD$ is itself reducible.
\end{proof}

\subsection{Families of rank 4 matrices}

We consider first the families of aCM bundles of type (\ref{4-4}),
We will see that their behavior is essentially the same as type (\ref{4-3}).
Recall that if $\GG$ is of type (\ref{4-4}) then $c_1(\GG)=0$ and $c_2(\GG)=1$,
while $c_1(\EE)=H$ and $c_2(\EE)=2$ if $\EE$ is of type (\ref{4-3}).

\begin{thm} \label{moduli-4-0}
Let $\GG$ and $\EE$ be indecomposable aCM bundles respectively of type (\ref{4-4}) and (\ref{4-3}).
Then we have:
\begin{enumerate}[i)]
\item \label{semistable-G} the bundle $\GG$ is strictly semistable unless it is an extension of the form:
	\begin{equation} \label{L-L-ext}
	0 \rr \OO_X(L) \rr \EE \rr \OO_X(-L) \rr 0 \qquad \mbox{with $L \in \LX$}
	\end{equation}
\item \label{stable-E} the bundle $\EE$ is either stable or an extension of the form:
	\begin{equation} \label{C-L-ext}
	0 \rr \OO_X(C) \rr \EE \rr \OO_X(H-C) \rr 0 \qquad \mbox{with $C \in \CX$}
	\end{equation}
\item \label{moduli-G} the moduli space $\MCMss(2;0,1)$ is isomorphic to $X$ via the map $\xi$;
\item \label{moduli-E} the correspondence $\phi: \GG \mapsto \ker(\ppp(\GG))(2)$ defines an isomorphism:
	\[ \MCMss(2;0,1) \rr \MCMst(2;H,2) 
	\]
\end{enumerate}
\end{thm}

\begin{proof}
The discussion in the proof of Proposition \ref{orientable-4} implies that a section
(unique up to a nonzero scalar) of $\GG$ vanishes along a single point $z$ of $X$,
unless $\GG$ is of the form (\ref{L-L-ext}). Let us work out the non-extension case. 
We have an exact sequence:
\begin{equation} \label{point-G}
0 \rr \OO_X \rr \GG \rr J_{z} \rr 0
\end{equation}

This exact sequence amounts to the Jordan-H\"older filtration of $\GG$.
This means that $\GG$ is semistable and in fact its $S$-equivalence class corresponds
to the class $[\OO_X \oplus J_{z}]$.
This gives the map $\xi : \MCMss(2;0,1) \rr X$.
On the other hand, any point $z\in X$ satisfies the Cayley-Bacharach property with
respect to $\OO_X(-H)$, and one sees easily that the extension provided by Theorem \ref{hartshorne-serre}
is an indecomposable aCM bundle of type (\ref{4-4}). So we have proved (\ref{semistable-G}) and (\ref{moduli-G}).

\medskip

Now let us turn to $\EE$. Combining (\ref{point-G}) with the minimal graded free resolution of $\GG$ we
can write down the following diagram (we omit zeroes all around):
\[
\xymatrix{
	\OO_X(1)\ar@{=}[d] \ar[r] & \GG(1) \ar[r] & J_{z}(1) \\
	\OO_X(1) \ar[r] & \OO_X(1) \oplus \OO_X^3 \ar[u]^-{\ppp(\GG)_{|X}} \ar[r] & \OO_X^3 \ar[u] \\
	& \EE^* \ar[u]^-{\ppp(\EE^*)^{\top}_{|X}} \ar@{=}[r] & \EE^* \ar[u]
	}
\]

Making use of the rightmost column of this diagram, we can prove (\ref{stable-E}).
Indeed a destabilizing subbundle $\KK^*$ of $\EE^*$ must have degree $0$ or $1$, with $\HH^0(\KK) \neq 0$.
In the former case $\EE$ is decomposable. In the latter we have $c_1(\KK) = L \in \LX$ and by a Chern
class computation we find the exact sequence (\ref{C-L-ext}).
Clearly, the correspondence $\phi$ takes an extension of the form (\ref{L-L-ext})
into one of the form (\ref{C-L-ext}). 
This means that it also takes $\MCMss(2;0,1)$ to $\MCMst(2;H,2)$. The map $\phi$ is obviously invertible.
\end{proof}

We take now into account bundles $\EE$ of type (\ref{4-1}).
Recall that in this case $c_1(\EE)=L+H$, $c_2(\EE)=2$, for some $L\in \LX$.
By the isomorphism (\ref{E-dual-F-2}) of 
Proposition \ref{chern-4-2}, this determines the behavior of bundles of
type (\ref{4-2}) as well. 

\begin{thm}
Let $\EE$ be of type (\ref{4-1}). Then $\EE$ is semistable.
The moduli space $\MCMss(2;L+H,2)$ is a smooth rational curve,
containing an open dense subset of stable bundles.
\end{thm}

\begin{proof}
The proof does not differ much from that of Proposition \ref{semistable-A}.
Introduce the same notation, and consider a destabilizing subbundle $\KK^*$.
We see immediately that $c_1(\KK)$ is an element $L'$ of $\LX$.
The aCM condition implies $L \neq L'$.
Rephrasing diagram (\ref{no-two-lines}) we discover that $\EE$ is decomposable if $L \cdot L'=1$.
On the other hand if $L \cdot L'=0$, we set $T=L+H-L'$, and check that $T$ lies in $\TX$.
Since $T \cdot L'=2$, we easily get that $\EE$ is decomposable applying Lemma \ref{its-extension}.

Now we define the open subset $\mathcal{H}_2^{\circ}$ of $\Hilb_2(X)$ by requiring that the subscheme $Z$
is contained in no divisor of degree $2$. Theorem \ref{hartshorne-serre} provides us with a map
$\zeta : \mathcal{H}_2^{\circ} \rr \FMCMst(2,L+H,2)$. 
Our discussion implies that this map is an isomorphism, whose inverse is $\xi$.
So the moduli space of these framed aCM bundles is a smooth rational variety of dimension $4$.
This proves our claim, since this space projects onto $\MCMst(2,L+H,2)$ with a $\p^3$ as generic fibre.
\end{proof}

\begin{prop}
Let $\EE$ be indecomposable of type (\ref{4-1}) and set $L = c_1(\EE)-H \in \LX$.
Then $\EE$ is strictly semistable if and only if it is an extension
of the form:
\begin{equation} \label{C-C-ext}
0 \rr \OO_X(C_1) \rr \EE \rr \OO_X(C_2) \rr 0
\end{equation}
for some $C_i \in \CX$, with $C_1 \cdot C_2 = 2$, $C_1 + C_2 = c_1(\EE)$.
\end{prop}

\begin{proof}
Take a general section of $\EE$ and consider its vanishing locus $Z$, which
consists of two distinct points of $X$.
A destabilizing subbundle $\KK$ must be effective of degree $2$,
and there must be a curve $\CC$ in $|\KK|$ containing $Z$.
By the discussion in the previous Theorem, we are reduced to the two cases
$c_1(\KK) = L_1 + L_2$, with $L_i \in \LX$, with $L_1 \cdot L_2 = 0$, or $c_1(\KK) = C_2 \in \CX$.

Let us consider the former. We have seen in the proof of the previous Theorem that $Z$ cannot be contained
in a line so we have $\len(Z \cap L_i)=1$, for $i=1,2$.
Writing down a diagram like (\ref{no-two-lines}), we see that $\EE$ is decomposable unless $L = L_i$,
for some $i$. But the Cayley-Bacharach property implies the contrary.

In the latter case, again a diagram similar to (\ref{no-two-lines}) implies that $\EE$ is decomposable
unless $C_2 \cdot L = 0$. In this case we set $C_1 = H+L-C_2$ and check that $C_1$ lies in $\CX$.
We obtain $C_1 \cdot C_2 = 2$. By Lemma \ref{its-extension} we obtain our statement if
$\CC$ is smooth. We leave it to the reader to verify the statement in case $\CC$ is reducible.
\end{proof}

\subsection{Rank 3 matrices}

The case of resolutions of type (\ref{3}) is summarized by the
following Proposition.

\begin{prop}
The isomorphism classes of bundles $\EE$ on $X$ admitting a
minimal resolution of the form (\ref{3}) are in one-to-one correspondence
with the $72$ elements of $\TX$ via the map $\EE \mapsto c_1(\EE)$.
Each bundle $\EE$ is stable and rigid. A general section of $\EE$ vanishes on a single point.
\end{prop}

\begin{proof}
Clearly, the set of bundles of the form (\ref{3}) is in one-to-one
correspondence with $\TX$ via the application:
\[ T \mapsto \EE := \ker(\ppp(\OO_X(T))_{|X})(2) \]

In turn this application agrees with the map $c_1$ defined above.
Since $\EE$ is globally generated, a general global section $s$ of $\EE$
vanishes on a single point $z$ of $X$, indeed $c_2(\EE)=c_1(T)^2=1$.
Considering the Koszul complex of the section $s$, and taking a
destabilizing line bundle $\KK$, we get $c_1(\KK) = L \in \LX$.
So $z$ lies in $L$.
Reasoning like in the proof of Proposition \ref{semistable-A}, we see that $\EE$
is decomposable if $T \cdot L \neq 0$, by $\HH^1(\OO_L(-z+(T-H)\cdot L))= 0$.
However if $T \cdot L=0$, then $T-L \in \CX$, and by Lemma \ref{its-extension}
we see easily that $\EE$ is decomposable.
The bundle $\EE$ is rigid since we have $\chi(\EEnd(\EE))=1$.
\end{proof}

\section{Counting Families} \label{numerology}

In this section we set up a few remarks in order to count the families we have encountered so far.
The next Theorem summarizes our results, where the next table has the following meaning.
In the first column {\em ref} we refer the reader to our previous results where we classify
these bundles according to their resolution and their Chern classes.
The column {\em stab} tells whether in each family we can find a strictly semistable ({\em ss}),
a stable ({\em st}) and a simple ({\em si}) bundle.
In the column {\em families} we write the dimension and the number of of each family.
The column {\em extensions} explains how many extensions of the form (\ref{the-extension}) there are
in each family. We denote $M = c_1(\MM)$, $N = c_1(\NN)$, $\hh^1 = \hh^1(\MM \ts \NN^*)$,
and {\em stab} indicates whether the relevant extension is (semi)stable or not.

\begin{thm} \label{all-numbers}
The families of rank 2 indecomposable aCM bundles on $X$ behave according to the following table.
\begin{tiny}
\begin{equation} \label{pallottoliere}
\begin{array}{c|c|c|c}
 &  \mbox{\em Stability} & \mbox{\em Families} & \mbox{\em Extensions} \\
\hline
\hline
\begin{array}{c}
\mbox{Ref.} \\
\hline
\hline
(\ref{6}.\ref{6-1}) \\
(\ref{6}.\ref{6-2}) \\
(\ref{6}.\ref{6-3}) \\
\hline
(\ref{5-2}.\ref{5-1-1}\,u) \\
(\ref{5-2}.\ref{5-1-1}) \\
(\ref{5-2}.\ref{5-1-2}) \\
\hline
(\ref{5-4}.\ref{5-2-1}\,u) \\
(\ref{5-4}.\ref{5-2-1}) \\
(\ref{5-4}.\ref{5-2-2}) \\
\hline
(\ref{4-1}) \\
(\ref{4-1}) \\
\hline
(\ref{4-3}\,u) \\
(\ref{4-3}) \\
(\ref{4-3}) \\
\hline
(\ref{4-2}) \\
(\ref{4-2}) \\
\hline
(\ref{4-4}\,u) \\
(\ref{4-4}) \\
\hline
(\ref{3}) \\
(\ref{3}) 
\end{array}
&
\begin{array}{c|c|c}
\mbox{\em ss} & \mbox{\em st} & \mbox{\em si} \\
\hline
\hline
\checkmark & \checkmark & \checkmark \\
\checkmark & \checkmark & \checkmark \\
\checkmark & \checkmark & \checkmark \\
\hline
 & & \checkmark \\
 & \checkmark & \checkmark \\
 & \checkmark & \checkmark \\
\hline
 &  & \checkmark \\
 & \checkmark & \checkmark \\
 & \checkmark & \checkmark \\
\hline
 \checkmark & \checkmark & \checkmark \\ 
 & \checkmark & \checkmark \\ 
\hline	
 & & \checkmark \\
 & \checkmark & \checkmark \\
 & \checkmark & \checkmark \\
\hline
 \checkmark & \checkmark & \checkmark \\ 
 \checkmark & \checkmark & \checkmark \\ 
\hline
 & &  \\ 
 \checkmark & & \\
\hline
 & \checkmark & \checkmark \\
 & \checkmark & \checkmark
\end{array}
& 
\begin{array}{c|c}
\mbox{\em num.} & \mbox{\em dim.} \\
\hline
\hline
1 & 5 \\
72 & 3 \\
270 & 1 \\
\hline
27 & 0 \\
27 & 2 \\
216 & 0 \\
\hline
27 & 0 \\
27 & 2 \\
216 & 0 \\
\hline
27 & 1 \\
27 & 1 \\
\hline
1 & 0 \\
1 & 2 \\
1 & 2 \\
\hline
27 & 1 \\
27 & 1 \\
\hline
1 & 0 \\
1 & 2 \\
\hline
72 & 0 \\
72 & 0
\end{array}
&
\begin{array}{c|c|c|c|c|c}
M & N & M \cdot N & \hh^1 & \mbox{\em num} & \mbox{\em stab} \\
\hline
\hline
T_1 & T_2 & 5 & 3 & 72 & \mbox{ss} \\
T_1 & T_2 & 4 & 2 & 20 & \mbox{ss} \\
T_1 & T_2 & 3 & 1 & 4 & \mbox{ss} \\
\hline
T & C & 3 & 1 & 16 & \mbox{u} \\
C & T & 3 & 2 & 16 & \mbox{st} \\
C & T & 2 & 1 & 5 & \mbox{st} \\
\hline
L & T-H & 1 & 1 & 16 & \mbox{u}  \\
T-H & L & 1 & 2 & 16 & \mbox{st}  \\
T-H & L & 0 & 1 & 5& \mbox{st}  \\
\hline
C_1 & C_2 & 2 & 1 & 10 & \mbox{ss}\\
L & T 		& 2 & 2 &16 & \mbox{ss}\\
\hline
C & L 				& 2 & 1 & 27 & \mbox{u} \\
T_1 - H & T_2 & 2 & 3 & 72 & \mbox{st} \\
L & C 				& 2 & 2 & 27 & \mbox{st} \\
\hline
L_1 & L_2 & 1 & 1 & 27 & \mbox{ss} \\
T-H & C & 1 & 2 & 16 & \mbox{st} \\
\hline
L & -L & 1 & 1 & 27 & \mbox{u} \\
-L & L & 1 & 2 & 27 & \mbox{ss} \\
\hline
L & C & 1 & 1 & 6 & \mbox{st} \\
T_1-H & T_2-H & 1 & 2 & 20 & \mbox{st}
\end{array}
\end{array}
\end{equation}
\end{tiny}
where $L,L_i \in \LX$, $C,C_i \in \CX$, $T,T_i \in \TX$.
\end{thm}

\begin{rmk}
It should be noted that bundles coming from an extension like (\ref{the-extension})
are strictly semistable whenever $\deg(\MM)=\deg(\NN)$. In this case, for each element of
$\HH^1(\MM \ts \NN^*)$ we obtain nonisomorphic extensions which however represents the same point
in $\MCMss(2;M+N,M\cdot N)$, so for instance for bundles of type (\ref{6}.\ref{6-1}) we get
$36$ semistable points in $\MCMss(2;2\,H,5)$.
By contrast if $\deg(\MM) < \deg(\NN)$ we obtain a projective space corresponding to $\p(\HH^1(\MM \ts \NN^*))$
sitting inside $\MCMss(2;M+N,M\cdot N)$.
\end{rmk}

\begin{proof}[Proof of \ref{all-numbers}].
We only have to enumerate the families and the extensions. 
According to the results of section \ref{moduli} all families containing a bundle $\EE$
are either irreducible
open dense subsets of (a union of components of) the moduli space $\Most(2;c_1(\EE),c_2(\EE))$ corresponding
to aCM sheaves, or a finite number of unstable bundles arising as extensions of aCM line bundles.
In any case it suffices to enumerate the relevant Chern classes.

It is thus straightforward to compute the number of all families, perhaps with the exception of
bundles of type (\ref{6}.\ref{6-3}) and (\ref{5-2}.\ref{5-1-2}).
For the first number, notice that $L \cdot C = 0 \Leftrightarrow L \cdot \rho(L)=1$,
and for each $L$, $\#\{L' \,|\, L' \cdot L = 1\} = 10$ so we get $27 \cdot 10 = 270$.
For the second: $\#\{L' \,|\, L' \cdot L = 0\} = 16$ and $27 \cdot 16/2 = 216$.

Turning to the number of extensions, we consider the finite maps:
\begin{scriptsize}
\begin{align}
\label{number-A2} & \{\{T_1,T_2\} \, | \, T_1 \cdot T_2 = 4 \} \xr{10:1} \TX && \{T_1,T_2\} \mapsto T_1+T_2-H \\
\label{number-A3} & \{\{T_1,T_2\} \, | \, T_1 \cdot T_2 = 3 \} \xr{4:1} \{(L,C) \, | \, L \cdot C = 0\}
	&& \{T_1,T_2\} \mapsto T_1+T_2-H \\ 
\label{number-B1} & \{(T,C) \, | \, T \cdot C = 3 \} \xr{16:1}  \CX && (T,C) \mapsto T+C-H \\ 
\label{number-B2} & \{(T,C) \, | \, T \cdot C = 2 \} \xr{5:1}  \{\{L_1,L_2\} \, | \, L_1 \cdot L_2 = 0\} && (T,C) \mapsto T+C-H \\ 
\label{number-D} & \{\{C_1,C_2\} \, | \, C_1 \cdot C_2 = 2 \} \xr{8:1}  \LX && \{C_1,C_2\} \mapsto C_1+C_2-H
\end{align}
\end{scriptsize}
where $L,L_i \in \LX$, $C,C_i \in \CX$, $T,T_i \in \TX$.
It is an easy but tedious exercise to check the cardinality of the fibres of the maps above;
the scrupulous reader may derived them from the tables in the Appendix \ref{app-2}.
The maps (\ref{number-A2}), (\ref{number-A3}), (\ref{number-B1}), (\ref{number-B2}), (\ref{number-D}) take care respectively
of cases (\ref{6}.\ref{6-1}), (\ref{6}.\ref{6-2}), (\ref{5-2}.\ref{5-1-1}), (\ref{5-2}.\ref{5-1-2}), (\ref{4-1}).
All the remaining cases can be obtained by duality.
\end{proof}

\appendix

\section{The blow up at six points of the projective plane} 

All the material contained in this appendix is well know and
we enclose it here for the reader's convenience.
We will actually use a tiny bit of the rich geometry coming into play when dealing
with cubic surfaces, such as Steiner triads, tritangent trios and so forth.
The interested reader can look, for instance, at the beautiful notes \cite{dolgachev:classical}.

\subsection{The $\mathbf{E}_6$ lattice}

Let $X$ be a smooth cubic surface.
The intersection product defines a lattice structure of signature $(1,6)$
on the group $\Pic(X) \simeq \Z^7$.
We write $\Z^{1,6}$ for the canonical $(1,6)$ lattice over the basis $(e_0,\ldots e_7)$, and set
$\kappa = 3\,e_0 - e_1 - \ldots - e_6$.
One defines the $\mathbf{E}_6$ lattice as $\mathbf{E}_6 = \kappa^{\bot} \subset \Z^{1,6}$.
A vector $v$ in $\mathbf{E}_6$ with resp. $(v,v)=-2$ is called a {\em root};
there are $72$ of them.
An {\em exceptional} vector is an element of $\Z^{1,6}$ with $(v,v)=-1$, $(v,\kappa)=-1$.
It is well-known that there are $72$ roots and $27$ exceptional vectors in $\Z^{1,6}$.
A {\em sixer} is a sextuple of mutually orthogonal exceptional vectors.
One checks that sixes are in one-to-one correspondence with roots,
by $(v_1,\ldots,v_6) \mapsto 1/3 \,(-2\,\kappa - v_1 - \ldots -v_6)$.

The {\em canonical root basis} is defined by $\alpha_0 = e_0 - e_1 - e_2 - e_3$,  $\alpha_i = e_i - e_{i+1}$,
and the (opposite of the) matrix of the bilinear form in this basis equals the Cartan matrix of the Dynkin diagram
of type $\mathbf{E}_6$. We consider the Weil group ${\sf W}(\mathbf{E}_6)$ acting on $\Z^{1,6}$, generated by
reflections along the hyperplanes $\alpha_i^{\bot}$. It has order $51840$.
The proof of the following Theorem can be found e.g. in \cite[Part II, Theorem 10.1.10]{dolgachev-verra:classical}.

\begin{thm} \label{weil}
The group ${\sf W}(\mathbf{E}_6)$ acts transitively on the sets of roots, exceptional vectors and sixes.
\end{thm}

A classical result says that $X$ is the
blow up at six points $P_1,\ldots P_6$ of the projective plane $\p^2$ over $\kk$, so 
let $\sigma : X \rr \p^2$ be the blow-down morphism.
Let $\ell$ be the pull-back by $\sigma$ of the class of a line
in $\p^2$; $b_1,\ldots,b_6$ be the exceptional divisors on $X$ associated to the points
$P_1,\ldots,P_6$.
The hyperplane divisor $H$ on $X$ is defined by $H =3\, \ell -\sum b_i$ and we have $\omega_X \simeq \OO_X(-H)$.

Factorizing $\sigma$ into $6$
blow-ups of single points $\sigma_1, \ldots, \sigma_6$ (i.e. ordering the points $P_1, \ldots, P_6$)
defines a {\em geometric marking} on $X$, i.e. an isometry of lattices
$\phi : \Pic(X) \rr \Z^{1,6}$, with $\phi(-H)=\kappa$. 
Any two geometric markings define a {\em Cremona isometry}, i.e. an isometry of $\Pic(X)$ preserving
the canonical class, and the group of Cremona isometries is isomorphic to ${\sf W}(\mathbf{E}_6)$.

\subsection{Lines, conics and twisted cubics} \label{app-1}

Recall the $27$ classes of lines and conics defined on $X$, after fixing the blowing--down morphism $\sigma$.
\begin{footnotesize}
\begin{align*}
 & L_i = b_i &&  L_{i,j} = \ell - b_i -b_j &&  L^j = 2\,\ell - \sum_{i \neq j} b_i   \\
 & C_i = \ell-b_i &&  C^{i,j} = 2\,\ell-\sum_{k\neq i,j} b_k && C^j = 3\,\ell - \sum_{i\neq j} b_i - 2\, b_j 
\end{align*}
\end{footnotesize}

Recall also the $72$ classes of twisted cubics defined on $X$:


\begin{scriptsize}
\begin{align*}
& T_0 = \ell \\
& T_{i,j,k} = 2\,\ell-b_i-b_j-b_k \\
& T_i^j=3\,\ell - \sum_{k \neq i,j} b_k -2\,b_i \\
& T^{i,j,k} = 4\ell - b_i -b_j - b_k -  2\sum_{l\neq i,j,k} b_{l} \\
& T^0 = 5\,\ell - 2\sum b_i
\end{align*}
\end{scriptsize}


We have the formulas:
\begin{align}
& L_i + C^i = L_{i,j} + C^{i,j} = L^i+C_i = H \\
& T_0 + T^0 = T_{i,j,k} + T^{i,j,k} = T_i^j + T^i_j = 2\,H
\end{align}

Fixing a geometric marking on $X$, 
it is easy to establish the one-to-one correspondences:
\begin{align*}
 \{\mbox{roots in $\Z^{1,6}$}\}  & \leftrightarrow  \{\mbox{classes of twisted cubics in $X$}\} \\
 \{\mbox{exceptional vectors in $\Z^{1,6}$}\} & \leftrightarrow  \{\mbox{lines in $X$}\} 
\end{align*}
where the first assignment sends a root $v$ the divisor class $H-\phi^{-1}(v)$.

\begin{rmk} \label{cremona}
Relabelling the classes considered so far, in such a way that the intersection form is preserved,
and fixing the canonical class, amounts to choosing a different root basis for $\Z^{1,6}$.
In other words, it corresponds to choosing a different geometric marking for $X$, i.e.
to a Cremona isometry of $\Pic(X)$.
Now recall that the group of Cremona isometries is ${\sf W}(\mathbf{E}_6)$, which acts
transitively in the set of lines, conics and twisted cubics in $X$.

Therefore, in any statement concerning an arbitrary pair $(D_1,D_2)$ of such divisor,
we are allowed to fix $D_1$ and check the statement for all $D_2$'s.
\end{rmk}

\subsection{Intersection numbers} \label{app-2}

Let $L,L' \in \LX$, $C \in \CX$ and $T \in \TX$ and let $D \in \{L,C,T \}$.
Taking a divisor class $D' \in \LX \cup \CX \cup \TX$, we subdivide the sets of $\LX$, $\CX$, $\TX$
according to the intersection number $D \cdot D'$.
In view of Remark \ref{cremona}, we will let:
$$L=L_1=b_1 \qquad C=C_1=\ell - b_1 \qquad T=T_0=\ell$$

In the following tables, {\em Int.} denotes the intersection number
$D \cdot D'$.
{\em Num.} (resp. {\em Tot.}) denotes the number of classes $D'$
having the given intersection against $D$ and a fixed coefficient for
$\ell$ (resp. regardlessly of the coefficient for $\ell$).
We consider first the case $D=L=L_1$.

\begin{tiny}
\begin{align*}
& \fbox{$L \cdot L'$} && 
\left\{
\begin{array}{c|c||c|c|c}
  \text{Int.} & \text{Tot. $L'$} & \text{Class $L'$} & \text{Indices} & \text{Num. $L'$}\\
  \hline
  \hline
  -1 & 1 & L_1 & & 1 \\
  \hline
    0 & 16 & L_i & i \neq 1 & 5 \\
    0 & 16 & L_{i,j} & 1<i<j & 10 \\
    0 & 16 & L^1 &  & 1 \\
  \hline
    1 & 10 & L_{1,i} & i\neq 1 & 5 \\
    1 & 10 & L^i & i\neq 1 & 5
\end{array} 
\right.\\
& \fbox{$L \cdot C'$} && 
\left\{
\begin{array}{c|c||c|c|c}
  \text{Int.} & \text{Tot. $C'$} & \text{Class $C'$} & \text{Indices} & \text{Num. $C'$}\\
  \hline
  \hline
    0 & 10 & C_i & i \neq 1 & 5 \\
    0 & 10 & C^{1,i} & i \neq 1 & 5 \\
  \hline
    1 & 16 & C_1 &  & 1 \\
    1 & 16 & C_{i,j} & 1<i<j & 10 \\
    1 & 16 & C^i & i \neq 1 & 5 \\
  \hline
    2 & 1 & C_{1} & i\neq 1 & 1
\end{array}
\right. \\
& \fbox{$L \cdot T'$} && 
\left\{
\begin{array}{c|c||c|c|c}
  \text{Int.} & \text{Tot. $T'$} & \text{Class $T'$} & \text{Indices} & \text{Num. $T'$}\\
  \hline
  \hline
    0 & 16 & T_0 &  & 1 \\
    0 & 16 & T_{i,j,k} & 1 <i<j<k & 10 \\
    0 & 16 & T^1_i & i \neq 1 & 5 \\
  \hline
    1 & 40 & T_{1,i,j} & 1<i<j & 10 \\
    1 & 40 & T^i_j & 1\neq i\neq j \neq 1 & 20 \\
    1 & 40 & T^{1,i,j} & 1<i<j  & 10 \\
  \hline
    2 & 16 & T^i_1 & i\neq 1 & 5 \\
    2 & 16 & T^{i,j,k} & 1<i<j<k & 10 \\
    2 & 16 & T^0 & & 1
\end{array}
\right.
\end{align*}
\end{tiny}

We consider then $D=C=C_1$.

\begin{scriptsize}
\begin{align*}
& \fbox{$C \cdot L'$} && 
\left\{
\begin{array}{c|c||c|c|c}
  \text{Int.} & \text{Tot. $L'$} & \text{Class $L'$} & \text{Indices} & \text{Num. $L'$}\\
  \hline
  \hline
    0 & 10 & L_i & i \neq 1 & 5 \\
    0 & 10 & L_{1,j} & i \neq 1 & 5  \\
  \hline
    1 & 16 & L_1 &  & 1 \\
    1 & 16 & L_{i,j} & 1<i<j & 10 \\
    1 & 16 & L^i & i\neq 1 & 5 \\
  \hline
  2 & 1 & L^1 & & 1 \\
\end{array} 
\right.\\
& \fbox{$C \cdot C'$} && 
\left\{
\begin{array}{c|c||c|c|c}
  \text{Int.} & \text{Tot. $C'$} & \text{Class $C'$} & \text{Indices} & \text{Num. $C'$}\\
  \hline
  \hline
    0 & 1  & C_1 &  & 1 \\
  \hline
    1 & 16 & C_i & i \neq 1 & 5 \\
    1 & 16 & C^{i,j} & 1<i<j & 10 \\
    1 & 16 & C^1 &  & 1 \\
  \hline
    2 & 1 & C^{1,i} & i\neq 1 & 5 \\
    2 & 1 & C^{i} & i\neq 1 & 5
\end{array}
\right. \\
& \fbox{$C \cdot T'$} && 
\left\{
\begin{array}{c|c||c|c|c}
  \text{Int.} & \text{Tot. $T'$} & \text{Class $T'$} & \text{Indices} & \text{Num. $T'$}\\
  \hline
  \hline
    1 & 16 & T_0 &  & 1 \\
    1 & 16 & T_{1,i,j} & 1 <i<j   & 10 \\
    1 & 16 & T_1^i & i \neq 1 & 5 \\
  \hline
    2 & 40 & T_{i,j,k} & 1<i<j<k & 10 \\
    2 & 40 & T^i_j & 1\neq i\neq j \neq 1 & 20 \\
    2 & 40 & T^{i,j,j} & 1<i<j<k  & 10 \\
  \hline
    3 & 16 & T^1_i & i\neq 1 & 5 \\
    3 & 16 & T^{1,i,j} & 1<i<j & 10 \\
    3 & 16 & T^0 & & 1
\end{array}
\right.
\end{align*}
\end{scriptsize}

Finally, we write the intersection numbers in case $D=T=T_0$.

\begin{scriptsize}
\begin{align*}
& \fbox{$T \cdot L'$} && 
\left\{
\begin{array}{c|c|c}
  \text{Int.} & \text{Tot. $L'$} & \text{Class $L'$} \\
  \hline
  \hline
    0 & 6 & L_i \\
  \hline
    1 & 15 & L_{i,j} \\
  \hline
    2 & 6 & L^i \\
\end{array} 
\right.\\
& \fbox{$T \cdot C'$} && 
\left\{
\begin{array}{c|c|c}
  \text{Int.} & \text{Tot. $C'$} & \text{Class $C'$} \\
  \hline
  \hline
    1 & 6  & C_i \\
  \hline
    2 & 15 & C_{i,j} \\
  \hline
    3 & 6 & C^i \\
\end{array}
\right. \\
& \fbox{$T \cdot T'$} && 
\left\{
\begin{array}{c|c|c}
  \text{Int.} & \text{Tot. $T'$} & \text{Class $T'$} \\
  \hline
  \hline
    1 & 1 & T_0 \\
  \hline
    2 & 20 & T_{i,j,k} \\
  \hline
    3 & 30 & T^i_j \\
  \hline
    4 & 20 & T^{i,j,k} \\
  \hline
    5 & 1 & T^0
\end{array}
\right.
\end{align*}
\end{scriptsize}

\def\cprime{$'$} \def\cprime{$'$} \def\cprime{$'$} \def\cprime{$'$}
  \def\cprime{$'$}
\providecommand{\bysame}{\leavevmode\hbox to3em{\hrulefill}\thinspace}
\providecommand{\MR}{\relax\ifhmode\unskip\space\fi MR }
\providecommand{\MRhref}[2]{%
  \href{http://www.ams.org/mathscinet-getitem?mr=#1}{#2}
}
\providecommand{\href}[2]{#2}

\end{document}